\documentclass[invmat, envcountsect]{svjour}

\usepackage{epsfig, xypic, eucal, times, amssymb, latexsym}
\renewcommand{\cal}{\mathcal}

\def\makeheadbox{{%
\hbox to0pt{\vbox{\baselineskip=10dd\hrule\hbox
to\hsize{\vrule\kern3pt\vbox{\kern3pt
\hbox{\bf This and other preprints availible at:}
\hbox{\hspace{1cm} http://www.math.uchicago.edu/{\textasciitilde}nathand/}
\kern3pt}\hfil\kern3pt\vrule}\hrule}%
\hss}}}

\makeatletter
\long\def\@caption#1[#2]#3{\par\addcontentsline{\csname
  ext@#1\endcsname}{#1}{\protect\numberline{\csname
  the#1\endcsname}{\ignorespaces #2}}\begingroup
    \@parboxrestore
    \@makecaption{\csname fnum@#1\endcsname}{\ignorespaces #3}\par
  \endgroup}

\def\firstcaption{\refstepcounter\@captype\@dblarg%
            {\@firstcaption\@captype}}

\def\secondcaption{\refstepcounter\@captype\@dblarg%
            {\@secondcaption\@captype}}

\long\def\@firstcaption#1[#2]#3{\addcontentsline{\csname
  ext@#1\endcsname}{#1}{\protect\numberline{\csname
  the#1\endcsname}{\ignorespaces #2}}\begingroup
    \@parboxrestore
    \vskip10pt
    \@maketwocaptions{\csname fnum@#1\endcsname}{\ignorespaces #3}%
    \ignorespaces\hspace{.073\textwidth}\hfil%
  \endgroup}

\long\def\@secondcaption#1[#2]#3{\addcontentsline{\csname
  ext@#1\endcsname}{#1}{\protect\numberline{\csname
  the#1\endcsname}{\ignorespaces #2}}\begingroup
    \@parboxrestore
    \@maketwocaptions{\csname fnum@#1\endcsname}{\ignorespaces #3}\par
  \endgroup}

\long\def\@maketwocaptions#1#2{%
   \parbox[t]{.46\textwidth}{{\bf #1.} #2}}

\def\fig@type{figure}

\long\def\@makecaption#1#2{\if\@captype\fig@type\vskip 10pt\fi
   \setbox\@tempboxa\hbox{{\bf #1.} #2}%
 \ifdim \wd\@tempboxa >\hsize {\bf #1.} #2\par
 \else \hbox to\hsize{\hfil\unhbox\@tempboxa\hfil}\fi
 \if\@captype\fig@type\else\vskip5.5pt\fi}

\makeatother

\spnewtheorem{Theorem}{Theorem}[section]{\bf}{\it}
\spnewtheorem{Conjecture}[Theorem]{Conjecture}{\bf}{\it}
\spnewtheorem{Lemma}[Theorem]{Lemma}{\bf}{\it}
\spnewtheorem{Corollary}[Theorem]{Corollary}{\bf}{\it}
\spnewtheorem{Proposition}[Theorem]{Proposition}{\bf}{\it}
\spnewtheorem{Claim}{Claim}{\it}{\rm}

\spnewtheorem*{Remark}{Remark}{\it}{\rm}
\spnewtheorem*{Remarks}{Remarks}{\it}{\rm}

\spnewtheorem{FullTheorem}{Theorem}{\bf}{\it}

\spnewtheorem{DegreeTheorem}{Theorem}{\bf}{\it}

\spnewtheorem{VolumeTheorem}{Volume Rigidity Theorem}{\bf}{\it}

\spnewtheorem{SecTwoLemma}{Lemma}{\bf}{\it}

\spnewtheorem{SecTwoProposition}{Proposition}{\bf}{\it}

\spnewtheorem{SecTwoTheorem}{Theorem}{\bf}{\it}

\spnewtheorem{VolSecLemma}{Lemma}[subsection]{\bf}{\it}
\spnewtheorem{VolSecCor}{Corollary}[subsection]{\bf}{\it}

\newcommand{\halmos}{\hfill $\Box$}

\newcommand{\Q}{{\Bbb Q}}
\newcommand{\R}{{\Bbb R}}
\newcommand{\C}{{\Bbb C}}
\newcommand{\Z}{{\Bbb Z}}
\newcommand{\N}{{\Bbb N}}
\renewcommand{\H}{{\Bbb H}}

\newcommand{\sC}{{\cal C}}

\newcommand{\bdry}{\partial}
\newcommand{\injects}{\hookrightarrow}
\newcommand{\surjects}{\to}
\newcommand{\cross}{\times}
\newcommand{\union}{\cup}

\newcommand{\maps}{\colon\thinspace}
\newcommand{\tr}{ {\mathrm{tr}} }

\newcommand{\SO}{ {\mbox{\rm SO}} }
\renewcommand{\O}{ {\mbox{\rm O}} }
\newcommand{\SL}[2]{\mbox{\rm SL}_{#1} #2}
\newcommand{\PSL}[2]{\mbox{\rm PSL}_{#1} #2 }
\newcommand{\CP}{{ {\C}\, 
    \mbox{\fontfamily{cmr}\fontshape{n}\selectfont P} }}
\newcommand{\mtext}[1]{\quad\mbox{#1}\quad}
\newcommand{\vol}{ {\mathrm{vol}}}
\newcommand{\Vol}{ {\mathrm{Vol}}}
\newcommand{\Area}{ {\mathrm{Area}}}
\newcommand{\df}{{\mbox{\rm \small df}}}
\newcommand{\smear}{{\mathrm{smear\,}}}
\newcommand{\Str}{{\mathrm{Str}}}
\newcommand{\Stab}{{\mathrm{Stab}}}
\newcommand{\fix}{{\mathrm{fix}}}
\newcommand{\diam}{{\mathrm{diam}}}
\newcommand{\interior}{{\mathrm{int}}}
\newcommand{\pair}[1]{\left\langle #1 \right\rangle}
\newcommand{\mysmallmatrix}[4]{\left(\begin{array}{cc}
          \scriptstyle #1 & \scriptstyle #2  \\ \scriptstyle #3 &
          \scriptstyle #4
          \end{array}\right)
          }

\begin{document}

\title{Cyclic surgery, degrees of maps of character curves, and volume
  rigidity for hyperbolic manifolds.}

\author{ Nathan M.~Dunfield \thanks{This work was partially supported
    by an NSF Graduate Fellowship.} }

\institute{ Department of  Mathematics, University of Chicago, 5734
  S.~University Ave., Chicago, IL 60637, USA. \\
  email: nathand@math.uchicago.edu}

\authorrunning{Nathan M.~Dunfield}
\titlerunning{Property P, degrees of maps of character curves, and
  volume rigidity }

\date{\today}
\maketitle

\section{Introduction}

This paper proves the following conjecture of Boyer and Zhang: If a
small hyperbolic knot in a homotopy sphere has a non-trivial cyclic
surgery slope $r$, then it has an incompressible surface with
non-integer boundary slope strictly between $r-1$ and $r +1$.  I state
this result below as Theorem~\ref{full_theorem}, after giving the
background needed to understand it.  Corollary~\ref{simpler_statement}
of Theorem~\ref{full_theorem} is that any small knot which has only
integer boundary slopes has Property~P.  The proof of
Theorem~\ref{full_theorem} also gives information about the diameter
of the set of boundary slopes of a hyperbolic knot.

The proof of Theorem~\ref{full_theorem} uses a new theorem about the
$\PSL2{\C}$ character variety of the exterior, $M$, of the knot.  This
result, which should be of independent interest, is given below as
Theorem~\ref{degree_theorem}.  It says that for certain components of
the character variety of $M$, the map on character varieties induced
by $\bdry M \injects M$ is a birational isomorphism onto its image.
The proof of Theorem~\ref{degree_theorem} depends on a fancy version
of Mostow rigidity due to Gromov, Thurston, and Goldman.  The
connection between Theorem~\ref{degree_theorem} and
Theorem~\ref{full_theorem} is the techniques introduced by Culler and
Shalen which connect the topology of $M$ with its $\PSL{2}{\C}$
character variety.

I will begin with the background needed for
Theorem~\ref{full_theorem}.  Let $K$ be a {\em knot} in a compact,
closed, 3-manifold $\Sigma$, that is, a tame embedding $S^1 \injects
\Sigma$.  The {\em exterior} of $K$, $M$, is $\Sigma$ minus an open
regular neighborhood of $K$.  So $M$ is a compact 3-manifold whose
boundary is a torus.  Suppose $\gamma$ is a simple closed curve in
$\bdry M$.  We can create a closed manifold $M_\gamma$ from $M$ by
taking a solid torus and gluing its boundary to $\bdry M$ in such a
way that $\gamma$ bounds a disc in the solid torus ($M_\gamma$ depends
only on the isotopy class of $\gamma$).  The new manifold $M_\gamma$
is called a {\em Dehn filling} of $M$ or a {\em Dehn surgery} on $K$.
Recently, many people have studied what kinds of manifolds arise when
you do this, especially in the case where $K$ is a knot in $S^3$ (see
the surveys \cite{GordonSurvey,Luecke95}).

From now on, let $\Sigma$ be a homotopy sphere.  An interesting
question is: Are there any non-trivial Dehn surgeries on $K$ such that
the resulting manifold is a homotopy sphere?  Or more generally, where
the resulting manifold has cyclic fundamental group?  A knot is said
to have {\em Property~P} if no non-trivial Dehn surgery on it yields a
homotopy sphere.  By the Cyclic Surgery Theorem \cite{CGLS}, for fixed
$K$ there is at most one $\gamma$ other than the meridian for which
$M_\gamma$ is a homotopy sphere.  Gordon and Luecke
\cite{GordonLuecke89b} have shown that non-trivial Dehn surgery on a
knot in $S^3$ never yields $S^3$, and so the Poincar\'e Conjecture
would imply that every knot has Property P. On the other hand, there
are plenty of examples of knots with non-trivial cyclic surgeries.
Fintushel and Stern \cite{FintushelStern1} discovered that the
$(-2,3,7)$ pretzel knot has two non-trivial surgeries where the
resulting manifold is a lens space (see also
\cite{BoyerMattmanZhang}).  More examples of knots with non-trivial
cyclic surgeries are given in \cite{Luecke95,BoyerZhang97}.  This
famed pretzel knot was also the first example found to have a
non-integer boundary slope (I will define what boundary slopes are in
a moment).  Boyer and Zhang conjectured that there was a general
connection between cyclic surgeries and non-integer boundary slopes.

A properly embedded surface $(F, \bdry F) \injects (M, \bdry M)$ will
be called {\em incompressible} if $\pi_1(F) \rightarrow \pi_1(M)$ is
injective and not a 2-sphere bounding a ball.  An incompressible
surface $F$ which is not boundary parallel is called {\em essential}.
Isotopy classes of oriented simple closed curves in the torus $\bdry
M$ are in bijective correspondence with primitive elements of
$H_1(\bdry M, \Z)$; unoriented isotopy classes with pairs of primitive
elements $\{ \gamma, -\gamma \}$.  If $(\alpha, \beta )$ is a basis
for $H_1(\bdry M, \Z)$, the {\em slope} of $\gamma = a \alpha + b
\beta$ with respect to this basis is $\frac{a}{b} \in \Q \union
\{\infty\}$.  Note that the slope of an unoriented isotopy class is
well defined, and gives a bijection between unoriented isotopy classes
of simple closed curves and $\Q \union \{\infty\}$.  If $\bdry F$ is
non-empty, it consists of disjoint simple closed curves in $\bdry M$.
These curves must all be parallel, and so correspond to the same pair
of primitive elements $\{ \gamma, -\gamma \}$ in $H_1(\bdry M, \Z)$.
These are called the {\em boundary classes} of $F$, and their slope
the {\em boundary slope\/} of $F$.  Hatcher \cite{Hatcher82} has shown
that the number of boundary slopes is always finite.  As $M$ is the
exterior of a knot in a homotopy sphere, fix a meridian-longitude
basis $(\mu, \lambda)$, for $H_1(\bdry M)$; all slopes will be with
respect to this basis.  A slope is integral if it is in $\Z$.  Many
knots have only integral boundary slopes, e.g.~all 2-bridge knots
\cite{HatcherThurston}.  If the Dehn filling $M_\gamma$ has cyclic
fundamental group, then the slope of $\gamma$ is called a {\em cyclic
  surgery slope}.  A manifold is {\em small} if it does not contain a
{\em closed}, non-boundary parallel, incompressible surface.  A knot
is small if its exterior is.  A knot is {\em hyperbolic} if the
interior of $M$ admits a complete hyperbolic metric of finite volume.
Finally, here is the theorem which was conjectured by Boyer and Zhang
in \cite{BoyerZhang97}, and improves their partial results in this
direction:
\begin{FullTheorem}
Suppose $K$ is a small hyperbolic knot in a homotopy sphere which has
a non-trivial cyclic surgery slope $r$.  Then there is an essential
surface in the exterior of $K$ whose boundary slope is non-integral
and lies in $(r - 1, r + 1)$.
\end{FullTheorem}
Note that the Cyclic Surgery Theorem \cite{CGLS} shows that a cyclic
slope for a knot in a homotopy sphere is always integral.  The
$(-2,3,7)$ pretzel satisfies the hypotheses of this theorem, which
explains why it has a non-integral boundary slope.
Theorem~\ref{full_theorem} remains true for some knot exteriors in
manifolds with non-trivial cyclic fundamental group, see
Section~\ref{pf_of_full_theorem}.  A corollary of the theorem is a
condition which implies that a knot has Property P:
\begin{Corollary}\label{simpler_statement}
  If a small knot in a homotopy sphere has only integral boundary
  slopes, then it has Property P.
\end{Corollary}
This follows as a small non-hyperbolic knot is a torus knot, and
surgeries on these have been classified \cite{Moser71}.

I should mention that there are small hyperbolic knots which have
non-integral boundary slopes but no cyclic surgeries.  Boyer and Zhang
gave examples in \cite{BoyerZhang97}, and one can also use the proof
of Proposition~2.2 of \cite{HatcherOertel} to give examples of small
Montesinos knots with non-integral boundary slopes of absolute value
less than 2.  So it is not possible to use
Corollary~\ref{simpler_statement} to show that all small knots have
Property~P.

The ideas of the proof of Theorem~\ref{full_theorem} can also be used to 
give information about the diameter, $d$, of the set of boundary slopes.
In \cite{CullerShalenDiameter}, Culler and Shalen showed that for any
knot the diameter $d \geq 2$.   For a hyperbolic knot, I show that if $d = 2$ then the
greatest and least slopes are not integral (see Section~\ref{norm} for a  
detailed statement).

Let me change course and state the theorem about character varieties
that underlies the proof of Theorem~\ref{full_theorem}, and which I hope
will be of use in other situations.  Let $M$ be a finite-volume
hyperbolic 3-manifold with one cusp.  Let $X(M)$ denote the
$\PSL{2}{\C}$ character variety of $M$.  Basically, $X(M)$ is the set
of representations of $\pi_1(M)$ into $\PSL{2}{\C}$ mod conjugacy, and
is an algebraic variety over $\C$ (for details, see
Section~\ref{character_varieties}).  Culler and Shalen introduced a way
of getting essential surfaces from $X(M)$, which has been very
useful in proving theorems about Dehn surgery.  I will postpone
explaining how this works until the next section, but this is what
connects the next theorem with the preceding.  Let $X_0$ denote an
irreducible component of $X(M)$ which contains the conjugacy class of
a discrete faithful representation.  Since $M$ has one cusp, the
(complex) dimension of $X_0$ is 1.  The inclusion $i \maps \bdry M
\injects M$ induces a map on character varieties $i^*\maps X_0 \to
X(\bdry M)$.  I will prove:
\begin{DegreeTheorem}
The map $i^*\maps X_0 \to X(\bdry M)$ is a birational isomorphism onto
its image.
\end{DegreeTheorem}
The conclusion of this theorem does {\em not}\/ always hold for other
components of $X(M)$.  The key to the proof of Theorem
\ref{degree_theorem} is:
\begin{VolumeTheorem}[Gromov-Thurston-Goldman]
Suppose $N$ is a compact, closed, hyperbolic 3-manifold.  If $\rho
\maps \pi_1 (N) \to \PSL{2}{\C}$ is a representation with $\vol(\rho)
= \vol(N)$, then $\rho$ is discrete and faithful.
\end{VolumeTheorem}
If $N$ is a closed manifold, the volume of a representation $\rho \maps
\pi_1(N) \to \PSL{2}{\C}$ is defined as follows.  Choose any smooth
equivariant map $f \maps \widetilde{N} \to \H^3$.  The form
$f^*(\Vol_{\H^3})$ descends to a form on $N$.  The volume $\vol(\rho)$
is the absolute value of the integral of this form over $N$.The volume
is independent of $f$ as any two such maps are equivariantly homotopic
(see Section~\ref{volume_of_rep}).  Goldman \cite{Goldman82}
noticed that one can prove Theorem~\ref{volume_rigidity} in
essentially the same way as the Strict Version of Gromov's Theorem
given by Thurston in his lecture notes \cite{ThurstonLectureNotes}.
As this section of Thurston's notes remains unpublished, I will include
a proof of Theorem~\ref{volume_rigidity}.
(Theorem~\ref{volume_rigidity} is now known to hold for all connected
semisimple Lie groups, where in the definition of volume $\H^3$ is
replaced by the appropriate symmetric space.  The cases not done in
\cite{Goldman82} are covered by \cite{Corlette88} and
\cite{Corlette92}).
 
The Volume Rigidity Theorem is connected to
Theorem~\ref{degree_theorem} by a theorem which is a reinterpretation
of some of the results of \cite{CCGLS}.  It shows, roughly, that the
volume of a representation $\rho$ depends only on the image of its
conjugacy class under the map $i^* \maps X(M) \to X(\bdry M)$ (for a
precise statement, see Section~\ref{volume_form}).

Let me end the introduction with an outline of the proof of
Theorem~\ref{full_theorem}.  The proof goes by contradiction, and so
suppose that $K$ has a cyclic surgery slope but no non-integer
boundary slopes.  A simple algebraic argument determines exactly what
the image of $i^* \maps X_0 \to X(\bdry M)$ is. Theorem
\ref{degree_theorem} says that $X_0$ is essentially the same as
$i^*(X_0)$, and so $X_0$ is now known.  Reading off from $X_0$
information about the number of boundary components of a certain
essential surface contained in $M$ leads to a contradiction.

The rest of this paper is organized as follows: In
Section~\ref{character_varieties} I will review the facts about character
varieties that I will need later, and also prove some needed lemmas.
Sections~\ref{pf_of_degree_theorem} and \ref{pf_of_full_theorem} are
devoted to the proofs of Theorems~\ref{degree_theorem} and
\ref{full_theorem} respectively.  Section \ref{norm} discusses the
proof of Theorem~\ref{full_theorem} in the context of the norm
introduced in \cite{CGLS}, and gives an application to the question of
the diameter of the set of boundary slopes.  Finally,
Section~\ref{pf_of_volume_theorem} gives a proof of the Volume Rigidity
Theorem.

I would like to thank Peter Shalen for his encouragement and innumerable
useful conversations, as well as telling me about
Proposition~\ref{number_of_boundary_components}.  I thank Andrew
Przeworski and Marc Culler for their comments on an earlier draft of
this paper.  I also thank the referee for very useful comments and
suggestions, and especially for pointing out an error in the original
statement of Theorem~\ref{diameter}.

\section{Character varieties}\label{character_varieties}

\subsection{Basics}

In this section, I will review facts about character varieties that will
be needed later.  The basic references are the first chapter of
\cite{CGLS} and \cite{CullerShalen83}, as well as the expository
article \cite{ShalenHandbook}.  The case of $\PSL{2}{\C}$-, as opposed
to $\SL{2}{\C}$-, character varieties, is treated in
\cite{BoyerZhangSeminorm}.

If $M$ is a topological space with finitely generated fundamental
group, denote by $R(M)$ the set of representations of $\pi_1(M)$ into
$\PSL{2}{\C}$.  This set has a natural structure as an affine complex
algebraic variety.  Now $\PSL{2}{\C}$ acts on $R(M)$ by conjugation of
representations.  Let $X(M)$ denote the quotient space (strictly
speaking, the algebro-geometric quotient), and $t\maps R(M) \surjects
X(M)$ the quotient map.  The affine variety $X(M)$ is called the
character variety of $R(M)$, as two representations in $R(M)$ map to
the same point of $X(M)$ if and only if they have the same character.
If two representations have the same character and one of them is
irreducible, then they are conjugate.  Moreover, for an irreducible
component $X$ of $X(M)$ which contains the character of an irreducible
representation, characters of reducible representations form a
subvariety of strictly smaller dimension.  So usually just think of
$X(M)$ as representations mod conjugacy.  A character is called
discrete, faithful, or whatever if all representations with that
character are discrete, faithful, or whatever.  For each $\gamma$ in
$\pi_1 (M)$ there is a regular function $f_\gamma$ on $X(M)$ such that
if $\chi$ in $X(M)$ is the character of a representation $\rho$ then
$f_\gamma(\chi) = \tr(\rho(\gamma))^2 - 4$.  This is well defined as
the trace of an element of $\PSL{2}{\C}$ is defined up to sign.  A map
between spaces $f \maps M \to N$ gives a map $f^*\maps X(N) \to X(M)$
via composition with $f_* \maps \pi_1(M) \to \pi_1(N)$.  In
particular, $i \maps \bdry M \to M$ induces a map $i^* \maps X(M) \to
X(\bdry M)$.

We can also consider $\SL{2}{\C}$- rather than $\PSL{2}{\C}$-
representations, and everything works the same way.  Let
$\tilde{R}(M)$ denote the variety of representations of $\pi_1(M)$
into $\SL{2}{\C}$, and $\tilde{X}(M)$ the associated character
variety.  In this case, there are also regular functions $\tr_\gamma$
on $\tilde{X}(M)$ defined by $\tr_\gamma(\chi) = \tr(\rho(\gamma))$.
The natural map $\tilde{X}(M) \to X(M)$ is finite to 1, though it may
not be onto, as not all representations into $\PSL{2}{\C}$ lift to
ones into $\SL{2}{\C}$.

When $M$ is a finite-volume hyperbolic manifold with one cusp, I will
denote by $X_0$ a component of $X(M)$ which contains the character of
a discrete faithful representation (note $X_0$ may not be unique, see
Section~\ref{df_reps}).   By \cite{Culler86}, a torsion
free discrete faithful representation lifts to $\SL{2}{\C}$ and there
is a component $\tilde{X}_0$ of $\tilde{X}(M)$ which covers $X_0$.
Said another way, every representation in $X_0$ lifts to one into
$\SL{2}(\C)$.

\subsection{Associated actions on trees}

Culler and Shalen discovered a way to construct essential
surfaces from $X(M)$ or $\tilde{X}(M)$.  In this and the next section,
I will explain their method for an irreducible curve $X$ in
$\tilde{X}(M)$ which contains the character of an irreducible
representation.  Let $Y$ be a smooth projective model for $X$ and $p
\maps Y \to X$ the associated birational isomorphism.  The finite set
of points $Y \setminus p^{-1}(X)$ are called ideal points.  Culler and
Shalen showed how to associate to each ideal point an action of
$\pi_1(M)$ on a simplicial tree as follows: Choose a curve $R \subset
\tilde{R}(M)$ so that the Zariski closure of $t(R)$ is $X$ (I will say
more about how to choose $R$ later in Lemma~\ref{choosing_R}).  If
$\C(X)$ denotes the field of rational functions on $X$, then $t$
induces an inclusion of fields $\C(X) \injects \C(R)$.  The point $y$
determines a valuation $v$ on $\C(X)$ where $v(f)$ is the order of $f$
at $y$ (zeros getting positive valuation).  There is a valuation $w$ on
$\C(R)$ which extends $v$ in the sense that there is a positive
integer $d$ so that $w(f) = d \cdot v(f)$ for all $f \in \C(X)$.
Associated to the field $\C(R)$ and the valuation $w$ is its
Bruhat-Tits tree, which I will denote $T_y$.  This is a simplicial tree
on which $\SL{2}{\C(R)}$ acts by isometeries.  The translation length
of an element $A$ of $\SL{2}{\C(R)}$ acting on $T_y$ is $\min(0, - 2
w(\tr(A)))$ (Proposition~II.3.15 of \cite{MorganShalen84}).  There is
a tautological representation $P \maps \pi_1(M)\to\SL{2}{\C[R]}$,
where $\C[R] \subset \C(R)$ is the ring of regular functions, defined
by $P(\gamma) = A$ where $A$ is the matrix of regular functions so
that at any representation $\rho \in R$ , $A(\rho) = \rho(\gamma)$.
Composing $P$ and the action of $\SL{2}{\C(R)}$ on $T_y$ gives the
promised action of $\pi_1(M)$ on a tree.

In order to prove Theorem~\ref{full_theorem} it will be necessary to
know the relationship between the valuations $v$ and $w$ exactly.
I will show:
\begin{SecTwoProposition}\label{associated_tree}
  For each ideal point $y$ in $Y$ there is a simplicial tree $T_y$ on
  which $\pi_1(M)$ acts so that the translation length of $\gamma \in
  \pi_1(M)$ acting on $T_y$ is twice the degree of pole of
  $\tr_\gamma$ at $y$.
\end{SecTwoProposition}

Because of the formula for translation length above, to prove
Proposition~\ref{associated_tree} it suffices to show we can choose
$R$ so that $d=1$ (you may wish to skip ahead until this becomes
crucial in the proof of Theorem~\ref{full_theorem}).  Let $Z$ be a
smooth projective model for $R$, and $\bar{t} \maps Z \to Y$ be the
map induced by $t \maps R \to X$.  The map $\bar{t}$ is a holomorphic
branched covering of Riemann surfaces. If there is a $z \in
\bar{t}^{-1}(y)$ such that $\bar{t}$ is not branched at $z$, take $w$
to be the valuation determined by order at $z$, and then $d = 1$.  Thus
Proposition~\ref{associated_tree} will be proved by:

\begin{SecTwoLemma}\label{choosing_R}
Let $y$ be a fixed ideal point of $X$.  Then there is a curve $R
\subset\tilde{R}(M)$ such that the induced map $\bar{t} \maps Z \to Y$
has the property that for any $z$ in $\bar{t}^{-1}(y)$, $\bar{t}$ is
not branched at $z$.
\end{SecTwoLemma}

\begin{proof}
By Corollary~1.4.5 of \cite{CullerShalen83} there is an $\alpha \in
\pi_1(M)$ so that $\tr_\alpha$ has a pole at $y$. Since $\tr_\alpha$
is non-constant there is an irreducible character $\chi$ in $X$ such
that $\tr_\alpha(\chi) \neq 2$.  Let $\rho$ be a representation with
character $\chi$.  By Proposition~1.5.1 of \cite{CullerShalen83} there is a
$\beta \in \pi_1(M)$ such that $\rho$ restricted to the subgroup
$\left\langle \alpha, \beta \right\rangle$ generated by $\alpha$ and
$\beta$ is irreducible, or equivalently $\tr_{[\alpha, \beta]}(\rho)
\neq 2$.  Considering the regular function $\tr_{[\alpha, \beta]}$, we
see that all but finitely many characters $\chi \in X$ are irreducible
when restricted to $\left\langle \alpha, \beta \right\rangle$.  There
are two cases, depending on whether $\tr_\beta$ is constant or not.
Suppose $\tr_\beta$ is non-constant.  Then there is a representation
$\rho_0$ whose character is in $X$ so that $\rho_0$ restricted to
$\left\langle \alpha, \beta \right\rangle$ is irreducible and
$\rho_0(\alpha)$ and $\rho_0(\beta)$ are hyperbolic elements of
$\SL{2}{\C}$.  Then $\rho_0$ restricted to $\left\langle \alpha,
\beta^n \alpha \beta^{-n} \right\rangle$ is irreducible for large $n$.
This is because the fixed point set of $\rho_0(\beta^n \alpha
\beta^{-n})$ is $\rho_0( \beta^n)(\fix(\alpha))$, which is disjoint
from $\fix(\alpha)$ for large $n$.  So there is a conjugate $\gamma$
of $\alpha$ such that $\rho_0$ restricted to $\left\langle \alpha,
\gamma \right\rangle$ is irreducible.  Again, all but finitely many
characters of $X$ are irreducible when restricted to $\left\langle
\alpha, \gamma \right\rangle$.  Let $V$ be the subvariety of
$\tilde{R}(M)$ consisting of representations $\rho$ such that
\begin{equation}\label{eqn_defining_V}
\rho(\alpha) = \left(\begin{array}{cc} a & 1 \\ 0 & 1/a
\end{array}\right), \mtext{and} \rho(\gamma) = \left(\begin{array}{cc}
a & 0\\ s  & 1/a \end{array}\right), 
\end{equation}
for some $a, s \in \C$, with $a \neq 0$.  Any representation $\rho$
which is irreducible when restricted to $\left\langle \alpha, \gamma
\right\rangle$ and where $\rho(\alpha)$ is not parabolic is conjugate
to exactly two representations in $V$.  For a conjugate of such a
$\rho$ lying in $V$, there are two choices for $a$, and once $a$ is
fixed, $s$ is determined by $\tr(\rho( \alpha \gamma))$.  Any two
conjugates of $\rho$ which agree on $\left\langle \alpha, \gamma
\right\rangle$ are actually equal, as the stabilizer under conjugation
of any irreducible representation is $\{ I, -I \}$.  Thus all but
finitely many points of $X$ are in $t(V)$, so we can choose an
irreducible curve $R \subset V$ so that $t(R)$ is dense in $X$.  The
map $t \maps R \to X$ is generically either 1-to-1 or 2-to-1.  In the
former case, we are done.  In the latter case $R = V$ and it is enough
to show that $\bar{t}^{-1}(y) \subset Z$ consists of two points.  Fix
$z$ in $\bar{t}^{-1}(y)$.  We can define a regular function $a \maps R
\to \C$ by Eqn.~(\ref{eqn_defining_V}).  Since $\tr_\alpha$ has a pole
at $y$, $a$ must have either a pole or a zero at any point of
$\bar{t}^{-1}(y)$.  Choose a sequence of representations $\rho_j$
converging to $y$.  From Eqn.~(\ref{eqn_defining_V}), we have
$\rho_j(\alpha) = \mysmallmatrix{ a_j }{ 1 }{ 0 }{1/a_j}$.  Since
generically a point of $X$ has two inverse images in $R$, for all but
finitely many $j$ there are representations $\psi_j \in R$ which are
conjugate to, but not equal to, $\rho_j$.  Hence $\psi_j(\alpha) =
\mysmallmatrix{1/a_j} {1}{0}{a_j}$ because as mentioned above,
conjugates which agree on $\left\langle \alpha, \gamma \right\rangle$
are equal.  Then the $\psi_j$ converge to a point $z'$ in
$\bar{t}^{-1}(y)$.  But $z$ is not $z'$ since if $a$ has a zero at $z$
then $a$ has a pole at $z'$ and vice versa.  So $\bar{t}^{-1}(y)$
consists of two points and we are done.

If $\tr_\beta$ is constant, fix $\lambda \in \C$ so that $\lambda +
1/\lambda = \tr_\beta$.  Define $V$ to be the subset of $R$ given by:
\[
\rho(\alpha) = \left(\begin{array}{cc} a & s \\ 0 & 1/a
\end{array}\right), \mtext{and} \rho(\beta) = \left(\begin{array}{cc}
\lambda & 0 \\ 1 & 1/\lambda\end{array}\right) 
\]
where $a, s \in \C$ and proceed as before. \halmos

\end{proof}

\subsection{Associated surfaces}\label{associated_surfaces}

There is a dual surface to any action of $\pi_1(M)$ on a tree $T$ as
follows.  Let $\tilde{M}$ be the universal cover of $M$.  Choose an
equivariant map $f \maps \tilde{M} \to T$ which is transverse to the
midpoints of the edges of $T$.  If $\tilde{S}$ is the inverse image
under $f$ of the midpoints of the edges of $T$, then $\tilde{S}$ is an
equivariant family of surfaces in $\tilde{M}$ which descends to a
surface $S$ in $M$.  In Section~1.3 of \cite{CGLS} it is shown how any
such $f$ can be modified so that $S$ becomes essential.  If $y$
is an ideal point, this construction gives an essential surface
dual to the action on $T_y$.  This surface, $S$, is said to be
associated to $y$.  Note that $S$ need not be connected or unique up
to isotopy.

Suppose $M$ has torus boundary, and $\bdry M$ is essential in
$M$.  Suppose $y$ is an ideal point where for some peripheral element
$\gamma \in \pi_1(\bdry M)$, $\tr_\gamma$ has a pole (there may be ideal
points where this does not happen).  In this case, there is exactly
one slope $\{\alpha, -\alpha\}$ such that $\tr_\alpha$ is finite at
$y$ (again, see \cite{CGLS} for details).  Then $\alpha$ fixes a point
of $T_y$.  If $S$ is an essential surface associated to $T_y$,
then some component of $S$ has non-empty boundary with boundary
classes $\{ \alpha, -\alpha \}$.  Let $\beta$ be such that $(\alpha,
\beta)$ is a basis for $\pi_1(\bdry M)$.  I will need:

\begin{SecTwoProposition}[Culler and
  Shalen]\label{number_of_boundary_components} The surface associated
  to an ideal point $y$ can be chosen so that the number of boundary
  components is equal to twice the order of pole of $\tr_\beta$ at
  $y$.
\end{SecTwoProposition}

\begin{proof}
This is essentially part of Section~5.6 of \cite{CCGLS}.  By
Proposition~\ref{associated_tree}, twice the order of pole of
$\tr_\beta$ at $y$ is the same as the translation length, $l(\beta)$,
of $\beta$ acting on $T_y$.  First, I claim $|\bdry S | \geq
l(\beta)$, where $|\bdry S |$ is the number of components of $\bdry
S$.  Let $p \maps \tilde{M} \to M$ be the covering map. Pick an arc
$a$ in $\tilde{M}$ which is a lift of a loop in $M$ representing
$\beta$ and which intersects $p^{-1}(S)$ in $|\bdry S|$ points.  Hence
the image of $a$ in $T_y$ intersects the midpoints of $|\bdry S|$
edges of $T_y$.  One of the endpoints of the image of $a$ in $T_y$ is
sent to the other under the action of $\beta$.  Therefore the
translation length of $\beta$ acting on $T_y$ is at most $|\bdry S|$.

Now I will produce a surface associated to $y$ with at most $l(\beta)$
boundary components.  Choose a connected component of $p^{-1}(\bdry
M)$, $C$, which we identify with $\R^2$ so that $\alpha$ acts on $C$
by unit translation in the first coordinate, and $\beta$ acts by unit
translation in the second coordinate.  The abelian subgroup
$\pi_1(\bdry M)$ leaves invariant a unique line $L$ in $T_y$, and acts
on $L$ via translations.  An element $\gamma = a \alpha + b \beta\in
\pi_1(\bdry M)$ acts on $L$ by translation by $b \cdot l(\beta)$.
Choose a map $f \maps C \to L$ which is equivariant under the action
of $\pi_1(\bdry M)$ as follows. Let $Y$ be the second coordinate axis.
First, project onto the second coordinate to get a map from $C$ to
$Y$.  Then compose this with a linear map from $Y$ to $L$ that expands
by a factor of $l(\beta)$.  There is a unique extension of $f$ to all
of $p^{-1}(\bdry M)$ which is equivariant under $\pi_1(M)$.  Since
$T_y$ is contractible, we can extend $f$ to an equivariant map of all
of $\tilde{M}$ to $T_y$.  The dual surface $S$ has $| \bdry S | =
l(\beta)$.  Changing $f$ so that $S$ becomes essential does not
increase the number of boundary components, so there is a surface
associated to $y$ with $|\bdry S| \leq l(\beta)$.  As we also have
$|\bdry S | \geq l(\beta)$, $S$ must have exactly $l(\beta)$ boundary
components. \halmos
\end{proof}

\subsection{Associated plane curves}\label{plane_curve}

The authors of \cite{CCGLS} introduced a plane curve 
associated to the character variety.  This gives a nice set of
coordinates for the computations in the proof of
Theorem~\ref{full_theorem}.  Let $M$ be a compact 3-manifold with torus
boundary.  Let $\Delta$ be the set of diagonal representations of
$\pi_1(\bdry M)$ into $\SL{2}{\C}$.  For any $\gamma \in \pi_1(\bdry
M)$ there is a well defined eigenvalue function $\xi_\gamma \maps \Delta
\to \C^*$ which takes a representation $\rho$ to the upper left hand
entry of $\rho(\gamma)$.  Fixing a basis $( \alpha, \beta )$ for
$\pi_1(\bdry M)$, the pair of eigenvalue functions $(\xi_\alpha,
\xi_\beta)$ gives coordinates on $\Delta$, and allows us to identify it
with $\C^*\times \C^*$.   Said another way,  we can identify a point
$(a, b) \in \C^*\times \C^*$ with the representation $\rho$ such that:
\[
\rho(\alpha) = \left(\begin{array}{cc} a & 0 \\ 0 & 1/a
\end{array}\right), \mtext{and} \rho(\gamma) =
\left(\begin{array}{cc}
b & 0\\ 0  & 1/b \end{array}\right).
\]
I will say that $(a, b) = (\xi_\alpha, \xi_\beta)$ are the
 eigenvalue coordinates corresponding to the basis $(\alpha, \beta)$.

There is a natural map $t \maps \Delta \to \tilde{X}(\bdry M)$
which is onto and generically 2-to-1.  If $Y$ is a one-dimensional
subvariety of $\tilde{X}(\bdry M)$, we can take the closure of
$t^{-1}(Y) \subset \Delta \subset \CP^2$ to get a plane curve $D$
which is a double cover of $Y$.  If $\tilde{X}$ is an irreducible
component of $\tilde{X}(M)$ such that $i^*(\tilde{X}) \subset
\tilde{X}(\bdry M)$ is one-dimensional, we can associate the plane
curve $D(\tilde{X}) \equiv t^{-1}(i^*(\tilde{X}))$ to
$\tilde{X}$.  The union of $D(\tilde{X})$ over all
components $\tilde{X}$ of $\tilde{X}(M)$ with
$i^*(\tilde{X})$ one-dimensional is called the associated plane
curve and denoted $D_M$.  Note for a component $\tilde{X}_0$ of
$\tilde{X}(M)$ which contains the character of a discrete faithful
representation, $i^*(\tilde{X}_0)$ is one-dimensional by
Proposition~1.1.1 of \cite{CGLS}.

There is an alternate construction of $D_M$ given in
\cite{CooperLong96} that we will need for the next section.  Let
$\tilde{R}_U(M)$ be the subset of $\tilde{R}(M)$ consisting of
representations whose restriction to the subgroup $\pi_1(\bdry M)$ is
upper-triangular.  Then there is a map $i^* \maps \tilde{R}_U(M)
\to \Delta$ which sends $\rho$ to the pair consisting of the upper
left hand entries of $\rho(\alpha)$ and 
$\rho(\beta)$.  It is not hard to see that the union of the
$i^*(\tilde{R})$ over all components $\tilde{R} \subset
\tilde{R}_U(M)$ with $i^*(\tilde{R})$ one-dimensional is
exactly $D_M$.

\subsection{Volume of a representation}\label{volume_of_rep}

Let $N$ be a closed 3-manifold.  The volume, $\vol(\rho)$, of a
representation $\rho \maps \pi_1(N) \to \PSL{2}{\C}$ is defined as
follows.  Choose any smooth equivariant map $f \maps \tilde{N} \to
\H^3$.  The form $f^*(\Vol_{\H^3})$ on $\tilde{N}$ descends to a form
on $N$.  The volume $\vol(\rho)$ is the absolute value of the integral
of this form over $N$. Since any two such maps are equivariantly
homotopic, the volume is independent of $f$.  If $f$ and $g$ are two
such maps one can use the straight line homotopy $H$ defined by $H(p,
t) = t f(p) + (1-t)g(p)$, where this linear combination is along the
geodesic joining $f(p)$ to $g(p)$.

The purpose of this section is to define the volume of a representation
for a compact 3-manifold $M$ whose boundary is a torus.  The
definition given in the closed case does not work here.  As $\bdry M$
is non-empty, the value of the integral used to define the volume
depends on the choice of equivariant map.

For convenience, I will assume throughout that the torus $\bdry M$ is
incompressible in $M$.  Let $\tilde{M}$ be the universal cover of $M$
with $\pi \maps \tilde{M} \to M$ the covering map.  Let $\bar{M}$
denote the quotient space of $\tilde{M}$ obtained by collapsing each
component plane of $\pi^{-1}(\bdry M)$ to a point.  The points of
$\bar{M}$ coming from the collapsed components of $\pi^{-1}(\bdry M)$
will be denoted $\bdry \bar{M}$; the set $\bar{M} \setminus \bdry
\bar{M}$ will be denoted $\interior(\bar{M})$.  The action of
$\pi_1(M)$ on $\tilde{M}$ induces an action on $\bar{M}$.  This action
is free on $\interior(\bar{M})$, but each point of $\bdry \bar{M}$ is
stabilized by some peripheral subgroup of $\pi_1(M)$.  Let
$\bar{\H}^3$ denote the union of $\H^3$ with the sphere at infinity
$S^2_\infty$.  The idea for defining the volume of a representation
$\rho \maps \pi_1(M) \to \PSL{2}{\C}$ is to consider equivariant maps
$f \maps \bar{M} \to \bar{\H}^3$ which send $\bdry \bar{M}$ into
$S^2_\infty$ and send $\interior(\bar{M})$ into $\H^3$, and then
proceed as in the closed case.

Before I can define the precise types of maps to be considered, I need
to define some notation.  Let $T$ be a torus and fix a product
structure $T \cross [0, \infty]$ on a neighborhood of $\bdry M$ such
that $\bdry M = T \cross \{\infty\}$ (the reason for this choice of
closed interval will become clear in a moment).  This gives an
equivariant product structure on $\pi^{-1}(T\cross [0,\infty])$ in
$\tilde{M}$.  This product structure induces an equivariant cone
structure on a neighborhood of $\bdry \bar{M}$ in $\bar{M}$.  If $v$
is in $\bdry \bar{M}$, I will denote the component of this
neighborhood containing $v$ by $N_v$, and the cone structure on $N_v$
by as $P_v \cross [0,\infty]$, where $P_v$ is a plane that covers $T
\cross \{0\}$ and $P_v \cross \{\infty\}$ is really just the cone
point $v$.  The peripheral subgroup of $\pi_1(M)$ which fixes $v$ will
be denoted $\Stab(v)$.  Note that $\Stab(v)$ preserves $N_v$.

Let $\rho \maps \pi_1(M) \to \PSL{2}{\C}$ be a representation.  A {\em
  pseudo-developing map\/} for $\rho$ is a smooth equivariant map $f
\maps \bar{M} \to \bar{\H}^3$ which sends $\bdry \bar{M}$ into
$S^2_\infty$ and $\interior(\bar{M})$ into $\H^3$ and which satisfies
the following additional condition.  For each $v$ in $\bdry \bar{M}$,
I require that $f$ maps each ray $\{p\} \cross [0,\infty]$ in the cone
neighborhood $N_v$ to a geodesic ray, and that $f$ parameterizes this
ray by arc length with respect to cone parameter in $[0, \infty]$ (by
continuity, this ray has endpoint $f(v)$ in $S^2_\infty$).  Note that
$f(v)$ must be a fixed point of $\rho(\Stab(v))$.  As an example, if
$\rho$ is the holonomy representation for a complete hyperbolic
structure on $\interior(M)$, then a developing map for this structure
extends to a pseudo-developing map of $\bar{M}$ by sending each $v$ in
$\bdry \bar{M}$ to the unique point in $S^2_\infty$ fixed by
$\rho(\Stab(v))$.  Another example is if we have a decomposition of
$M$ into ideal tetrahedra; an equivariant map which is piecewise
straight with respect to this triangulation is a pseudo-developing
map.  Every $\rho$ has a pseudo-developing map.  To construct one,
pick a $v$ in $\bdry \bar{M}$ and define $f(v)$ to be some point fixed
by $\rho(\Stab(v))$.  Extend $f$ to $N_v$ by picking any $\Stab(v)$
equivariant map of $P\cross \{0\}$ into $\H^3$ and then extending 
geodesically along the cone structure of $N_v$.  There is a unique
equivariant extension of $f$ to the cone neighborhood of $\bdry
\bar{M}$.  As $\H^3$ is contractible, $f$ extends to a
pseudo-developing map for $\rho$.

I will define the volume of a pseudo-developing map $f$ in the same
way as in the closed case.  The pull back $f^*( \Vol_{\H^3})$ descends
to a form on $M$.  The absolute value of the integral of this form
over $M$ is defined to be $\vol(f)$.  To see that this integral is
well defined, pick a $v \in \bdry \bar{M}$ and choose a fundamental
domain $F = D \cross [0,\infty]$ for the action of $\Stab(v)$ on
$N_v$.  We need that the integral of $f^*( \Vol_{\H^3})$ over $F$ is
defined and finite.  Pick a horoball in $\H^3$ centered at $f(v)$
which contains $f(F)$.  Project $f(D \cross \{0\})$ out to the
corresponding horosphere $S$ along geodesic rays with limit $f(v)$.
Call this projection $\theta$.  Let $\Area_S$ be the area form on $S$.
Let $A = \int_{D \cross \{0\}} \left| \theta^* \Area_S \right|$, the
total unoriented area of $\theta(D\cross \{0\})$.  As $f$ restricted
to $\{p\} \cross [0,\infty]$ is a geodesic parameterized by arc
length, the same computation as computing the volume of a cusp shows
that the integral of $\left| f^*(\Vol_{\H^3}) \right|$ over $F$ is
bounded by $A/2$.  Thus the integral of the absolute value of
$f^*(\Vol_\H^3)$ over $M$ is finite, and so $\vol(M)$ is well defined.
Note that for the two examples of pseudo-developing maps given in the
last paragraph, the volume of such a map is the volume that you would
expect.  It would also be possible to define the volume of $f$ without
taking the absolute value of the given integral.  All the lemmas in
this section remain true with this altered definition; this fact will
be needed in Section~\ref{volume_form}.

I would now like to say that the volume of a pseudo-developing map is
independent of the choice of map using the same argument as in the
closed case.  However, there is a problem.  Suppose $f$ is a
pseudo-developing map for a representation $\rho$ such that
$\rho(\pi_1(\bdry M))$ has exactly two distinct fixed points.
Consider a $v \in \bdry \bar{M}$.  Let $p$ and $q$ in $S^2_\infty$ be
the two fixed points of $\rho(\Stab(v))$.  Then $f(v)$ is either $p$
or $q$, say $p$.  We can construct a pseudo-developing map $g$ for
$\rho$ with $g(v) = q$.  Then $f$ and $g$ are not homotopic through
equivariant maps from $\bar{M}$ to $\bar{\H}^3$ which send $\bdry
\bar{M}$ into $S^2_\infty$.  I see no a priori way to compare the
volumes of these two maps.  However, we have:
\begin{VolSecLemma} 
If $f$ and $g$ are two pseudo-developing maps for $\rho$ which
 agree on $\bdry \bar{M}$, then $\vol(f) = \vol(g)$.
\end{VolSecLemma}
\begin{proof}
Construct a homotopy $H \maps \bar{M} \cross [0,1] \to \bar{\H}^3$
between $h$ and $g$ through pseudo-developing maps as follows.  Let
$v$ be in $\bdry \bar{M}$. Consider the cone $N_v = P_v \cross
[0,\infty]$.  Begin to construct $H$ by setting $H(v, t) = f(v) =
g(v)$ for all $t$.  Extend $H$ over $P_v \cross \{0\}$ by any
$\Stab(v)$ equivariant homotopy between the restrictions of $f$ and
$g$ to $P_v \cross \{0\}$.  As before, cone along geodesic rays ending
at $f(v)$ to extend $H$ over $N_v$.  We can now extend $H$ to the
desired homotopy.

To see that $\vol(f) = \vol(g)$ consider $M_t = M \setminus T \cross
(t, \infty]$, where $T\cross [0, \infty]$ is our collar on $\bdry M$.
For $t \in [0,\infty]$, define $V_t(f) = \int_{M_t} f^*(\Vol_{\H^3})$
and let $V_t(g)$ be the corresponding quantity for $g$.  Note
$|V_f(\infty)| = \vol(f)$ and similarly for $g$.  The form
$\Vol_{\H^3}$ is closed, so by Stokes Theorem:
\[
 0 = \int_{M_t \cross [0,1]} dH^*(\Vol_{\H^3}) = V_f(t) - V_g(t) +
 \int_{ (\bdry M_t) \cross [0,1]} H^*(\Vol_{\H^3})
\]
Because $H$ was constructed by geodesically coning over the
restriction of $H$ to $P_v \cross \{0\}$, the integral $\int_{ (\bdry
  M_t) \cross [0,1]} H^*(\Vol_{\H^3})$ goes to zero exponentially with
$t$.  Therefore $\vol(f) = \vol(g)$.  \halmos
\end{proof}

If $\rho$ is a representation such that the volume of any two
pseudo-developing maps agree, I will say that the volume of $\rho$ is
defined and set $\vol(\rho) = \vol(f)$ for any pseudo-developing map
$f$.  If $\rho$ is a representation where $\rho(\pi_1(\bdry M))$
contains a non-trivial parabolic, then $\rho(\pi_1(\bdry M))$ has a
unique fixed point in $S^2_\infty$.  By the lemma, any two
pseudo-developing maps for $\rho$ have the same volume, and so the
volume is defined for such $\rho$.   I will show:
\begin{VolSecLemma}
If $\rho$ is a representation whose character lies in an irreducible
component of $X(M)$ which contains the character of a discrete
faithful representation then $\vol(\rho)$ is defined.  
\end{VolSecLemma}
which follows immediately from:
\begin{VolSecLemma}
Suppose $\rho_t$, with $t \in [0,1]$ is a smooth one parameter family
of representations of $\pi_1(M)$ into $\PSL{2}{\C}$.  If the volume of
$\rho_0$ is defined then so is the volume of $\rho_1$.
\end{VolSecLemma}

\begin{proof}
I will need the results of Section 4.5 of \cite{CCGLS}.  There, they
construct a particularly nice form of pseudo-developing map as
follows.  Let $N$ denote $M$ with $\bdry M$ collapsed to a point.
Decompose $N$ as a simplicial complex so that each simplex has at most
one vertex at the collapsed $\bdry M$.  This induces an equivariant
simplicial decomposition of $\bar{M}$.  Pick a preferred vertex
$v_\infty$ of $\bdry {\bar M}$ and a fixed point $p$ of
$\rho(\Stab(v_\infty))$.  For a representation $\rho$ we can construct
a pseudo-developing map $f$ for $\rho$ so that $f$ is piecewise
straight on the simplices of $\bar{M}$ and which sends $v_\infty$ to
$p$.  Simply pick an equivariant map of the zero skeleton of
$\bar{M}$ sending $v_\infty$ to $p$ and extend linearly along the
simplices.  This gives a pseudo-developing map because each simplex
has at most one vertex in $\bdry \bar{M}$ and so the only points sent to
$S^2_\infty$ are those in $\bdry \bar{M}$.

Pick two pseudo-developing maps $f_1$ and $g_1$ for $\rho_1$.  We can
build a smooth family of maps $f_t$, for $t \in [0,1]$ where $f_t$ is a
pseudo-developing map for $\rho_t$ of the special type just discussed.
Let $g_t$ be a similar family of maps for $g_1$.  In Section 4.5 of
\cite{CCGLS} it is shown that work of Hodgson \cite{Hodgson86} implies
that the derivative of $\vol(f_t)$ (or $\vol(g_t)$) depends only on
the restriction of $\rho_t$ to $\pi_1(\bdry M)$.  This implies that
the derivatives of $\vol(f_t)$ and $\vol(g_t)$ are the same.  As
$\vol(f_0) = \vol(g_0) = \vol(\rho_0)$, we must have $\vol(f_1) =
\vol(g_1)$.  So $\vol(\rho_1)$ is well defined. \halmos
\end{proof}

For the proof of Theorem~\ref{degree_theorem}, I will need:

\begin{VolSecLemma}\label{factor_lemma}
Suppose $\rho$ is a representation of $\pi_1(M)$ which factors through
the fundamental group of a Dehn filling $M_\gamma$ of $M$.  Then the
volume of $\rho$ with respect to $M$ is defined, and is equal to the volume
of $\rho$ with respect to the closed manifold $M_\gamma$.
\end{VolSecLemma}

\begin{proof}
Let $C$ be the solid torus added to $M$ to make $M_\gamma$.  Let
$\tilde{M}_\gamma$ be the universal cover of $M_\gamma$, with $\phi$
the covering map.  Pick an equivariant map $f$ from $\tilde{M}_\gamma$
to $\H^3$.  Let $V = \tilde{M}_\gamma \setminus
\phi^{-1}(\interior(C))$.  Let $W$ be $\bar{M}$ minus the open cone
neighborhood of $\bdry{\bar{M}}$.  Adjusting collars, I will view $V$
as a quotient of $W$.  The restriction of $f$ to $V$ is a
$\pi_1(M_\gamma)$-equivariant map which induces a
$\pi_1(M)$-equivariant map $F \maps W \to \H^3$.  Extend $F$ over
$\bdry\bar{M}$ by any $\pi_1(M)$-equivariant map.  By coning along
geodesics we can extend $F$ over the cone neighborhood of
$\bdry\bar{M}$ to a pseudo developing map for $\rho$.

To compare volumes, choose $f \maps \tilde{M}_\gamma \to \H^3$ so that
$f(\phi^{-1}(C))$ is one-dimensional (this is possible because there
is map from $M_\gamma$ to itself which is the identity outside a
neighborhood of $C$ and which collapses $C$ to a core curve
in $C$).  Thus the volume of $\rho$ with respect to
$M_\gamma$ is the integral of $f^*(\Vol_{\H^3})$ over $M_\gamma
\setminus C$.  The image of the cone neighborhood of $\bdry\bar{M}$
under $F$ is at most two-dimensional, and so $\vol(f)$ is the integral
of $F^*(\Vol_{\H^3})$ over $M$ minus the collar on $\bdry M$.  Because
$F$ is a lift of $f$ these two integrals have the same value.  Thus
$\vol(f)$ is equal to the volume of $\rho$ with respect to $M_\gamma$.

Since the restriction of $F$ to $\bdry\bar{M}$ was an arbitrary
equivariant map, the volume of $\rho$ with respect to $M$ is defined
and is equal to the volume of $\rho$ with respect to $M_\gamma$.
\halmos
\end{proof}

\subsection{Volume form}\label{volume_form}

In this section $M$ will be a hyperbolic 3-manifold with one cusp.  In
Section~\ref{volume_of_rep}, I discussed the volume of a
representation of $\pi_1(M)$ into $\PSL{2}{\C}$.  As volume is
invariant under conjugation, there is a map $\vol$ from the
irreducible characters of $X(M)$ to $\R^+$ where $\vol(\chi)$ is the
volume of any representation with character $\chi$.  The next
proposition is a reinterpretation of the results of Sections 4.4-4.5
of \cite{CCGLS}, and will be one of the keys to the proof of
Theorem~\ref{degree_theorem}.  A normalization of a curve $X$ is a
smooth curve $Y$ together with a regular birational map $f \maps Y \to
X$.  A normalization of $X$ can be constructed by taking a smooth
projective model $Y'$ with birational isomorphism $f\maps Y' \to X$
and letting $Y = f^{-1}(X)$.

\begin{SecTwoTheorem} \label{factor_theorem}
Let $X_0$ be an irreducible component of $X(M)$ which contains the
character of a discrete faithful representation.  Let $Y$ be a
normalization of $i^*(X_0)$ where $i^* \maps X(M) \to X(\bdry M)$ is
the map induced by $i \maps \bdry M \injects M$.  Then the map $\vol$
on the irreducible characters of $X_0$ factors through a map $Y \to
\R^+$.
\end{SecTwoTheorem}

That is, there is a map $\vol \maps Y \to \R^+$ such that the diagram
\[
\xymatrix{ & & Y  \ar[dll]_\vol \cr
               {\R^+} & X_0 \ar[l]^\vol \ar[r]_{i^*} & i^*(X_0) \ar[u]_f
}
\] commutes.  In fact, $\vol$ will be the
absolute value of a generically smooth function from $Y$ to  $\R$.
 
\begin{proof}
I will use the notation of Section~\ref{plane_curve}.  
In Section~4 of \cite{CCGLS}, the authors define a real-valued
differential form $\eta$ on $\Delta$ by
\[
\eta = \log |a|\,d\arg(b) - \log |b|\,d\arg(a).
\]
which measures the change in volume in the following sense.  Let
$\tilde{R}_0$ denote the subset of $\tilde{R}_U(M)$ which maps to
$X_0$.  In the language of Section~\ref{volume_of_rep}, we can define
a map $V \maps \tilde{R}_0 \to \R$ by setting $V(\rho)$ to be the integral
of $f^*(\Vol_{\H^3})$ over $M$ where $f$ is any pseudo-developing map
for $\rho$.  Thus $V(\rho) = \pm \vol(\rho)$, and the results of
Section~\ref{volume_of_rep} show that $V$ is well defined.  Work of
Hodgson, see Section~4.5 of \cite{CCGLS}, shows that $V$ is smooth and
that $dV$ is the pull-back of $-\frac{1}{2}\eta$ along $i^* \maps
\tilde{R}_0 \to D_M$.  By either Section 4.4 or 4.5 of \cite{CCGLS}, the
form $\eta$ is exact on $D_M$, in the sense that it is exact (where
defined) on any smooth projective model of $D_M$.  Let $D'_M$ be
subset of $D_M$ on which $\eta$ is defined.  Therefore, $V$ factors
through a map from a normalization of $D'_M$ to $\R$.  More precisely,
let $\bar{D}_M$ be a normalization of $D'_M$, and $f\maps D'_M\to
\bar{D}_M$ a birational isomorphism.  Then there is a smooth function
$v \maps \bar{D}_M \to \R$ such that $v \circ f \circ i^*= V$.  Note
that the form $\eta$ is invariant under the transformations which
quotient $\Delta$ down to $X(\bdry M)$ (these transformations are $(a,
b) \mapsto (1/a, 1/b)$, $(a, b) \mapsto (-a , b)$, etc.).  Therefore
$v$ descends to a function $v'$ on a normalization of $t(D'_M)$.  As
any character in $X_0$ lifts to one in $\SL{2}{\C}$ and $\eta$ is
defined on all of $\Delta$, the curve $i^*(X_0)$ is contained in
$t(D'_M)$.  So the function $v'$ is defined on a normalization of
$i^*(X_0)$.  If $\chi \in X_0$ is the character of a representation
$\rho \in \tilde{R}_0$, we have $|v'(i^*(\chi))| = |V(\rho)| =
\vol(\rho) = \vol(\chi)$.  So $|v'|$ is the required function.
\halmos
\end{proof}

\subsection{Discrete faithful representations}\label{df_reps}

It is important to remember that if $\rho$ and $\rho' \maps \pi_1(M)
\injects \PSL{2}{\C}$ are holonomy representations of a one-cusped
finite-volume hyperbolic 3-manifold $M$, then $\rho$ and $\rho'$ need
{\em not} be conjugate in $\PSL{2}{\C}$.  They are conjugate in
$\O(3,1)$ by Mostow rigidity, and it's not difficult to see that
$\rho'$ is conjugate in $\PSL{2}{\C}$ to either $\rho$ or to the
complex conjugate $\bar{\rho}$ of $\rho$.  By complex conjugate, I
mean $\bar{\rho}(\gamma)$ is the matrix whose entries are the complex
conjugates of those of $\rho(\gamma)$, for all $\gamma \in \pi_1(M)$.
I will need the following little lemma, which the reader may wish to
skip until it is used in the proof of Theorem \ref{degree_theorem}.
It is just an application of Lemma 6.1 of \cite{CooperLong96}, and
describes the behavior of the map $i^* \maps X(M) \to X(\bdry M)$
near the two discrete faithful characters.  Let $p$ be the point of
$X(\bdry M)$ where the trace of any element of $\pi_1(\bdry M)$ is
$\pm 2.$ If $\chi$ is a discrete faithful character, then $i^*(\chi) =
p$.  The lemma is:

\begin{SecTwoLemma} \label{distinct_images_of_d.f._reps}
  Let $M$ be a finite-volume hyperbolic 3-manifold with one cusp.  
  Then the two discrete faithful characters have neighborhoods in 
  $X(M)$ whose images under $i^*$ are distinct branches of $i^*(X_0)$ 
  through $p$.
\end{SecTwoLemma}

\begin{proof} 
  Let $(\alpha, \beta)$ be a basis of $\pi_1(\bdry M)$.  If
  $\rho_0$ is a discrete faithful representation in $\PSL{2}{\C}$ we can
  conjugate it by an element of $\PSL{2}{\C}$ so that
\begin{equation}
  \rho_0(\alpha) = \pm \left(\begin{array}{cc} 1 & 1 \\ 0 & 1
\end{array}\right) \mtext{and} \rho_0(\beta) = \pm \left(\begin{array}{cc}
1 & a \\ 0 & 1 \end{array}\right) \label{distinct_images_eqn}.
\end{equation}
The cusp shape of $M$ with respect to this basis is $a$
(Eqn.~(\ref{distinct_images_eqn}) uniquely determines $a$).  The cusp
shape $a$ is always non-real, and, by changing the basis of
$\pi_1(\bdry M)$ if necessary, we can assume $a$ is not pure
imaginary.  Let $\chi_0 \in X(M)$ be the character of $\rho_0$.  As in
Lemma 6.1 of \cite{CooperLong96} it is not hard to show that
$\lim_{\chi \to \chi_0} f_\beta / f_\alpha = a^2$ (the key is to note
that $\alpha$ and $\beta$ commute).  From
Eqn.~(\ref{distinct_images_eqn}) we see that the cusp shape of
$\bar{\rho}_0$ is $\bar{a}$.  So, $\lim_{\chi \to \bar{\chi}_0}
f_\beta / f_\alpha = \bar{a}^2 \neq a^2$.  Since $f_\alpha$ and
$f_\beta$ depend only on the image of a character in $X(\bdry M)$,
neighborhoods of $\chi_0$ and $\bar{\chi}_0$ must go to distinct
branches of $i^*(X(M))$ through $p$.  \halmos
\end{proof}

\section{ Degrees of maps of character curves}
\label{pf_of_degree_theorem}

\begin{Theorem}\label{degree_theorem}
  Let $M$ be a finite-volume hyperbolic 3-manifold with one cusp.  Let
  $X_0$ be a component of the $\PSL{2}{\C}$ character variety of $M$
  which contains the character of a discrete faithful representation.
  The inclusion $i \maps \bdry M \injects M$ induces a regular map
  $i^* \maps X_0 \to X(\bdry M)$.  This map is a birational
  isomorphism onto its image.
\end{Theorem}

\begin{proof}

Fix a discrete faithful character $\chi^\df$ in $X_0$.  By
Proposition~1.1.1 of \cite{CGLS}, $X_0$ has complex dimension 1, and
$i^*$ is non-constant, so $i^*(X_0)$ also has dimension 1.  As $X_0$
is irreducible, so is $i^*(X_0)$.  The map $i^* \maps X_0 \to
i^*(X_0)$ is a regular map of irreducible algebraic curves, and so has
a degree which is the number of points in $(i^*)^{-1} (p)$ for generic
$p \in i^*(X_0)$ (here, generic means except for a finite number of
points of $i^*(X_0)$).  A degree-1 map is always a birational
isomorphism.  Thus it suffices to show that there are infinitely many
points $p$ in $i^*(X_0)$ where $(i^*)^{-1} (p)$ consists of a single
point.

We construct $p_j$, $j \in \N$, so that $(i^*)^{-1} (p_j)$ consists of
a single point as follows.  By Thurston's Hyperbolic Dehn Surgery
Theorem, all but finitely many Dehn fillings of $M$ are hyperbolic
(see \cite{ThurstonLectureNotes} and \cite{NeumannZagier}, or
\cite{BenedettiPetronio}).  Choose an infinite sequence of distinct
Dehn fillings so that the resulting manifolds $M_{\gamma_1},
M_{\gamma_2}, \ldots $ are all hyperbolic.  Moreover, Thurston's
theorem says that we can choose the $\gamma_j$'s so that the
hyperbolic structures of the $M_{\gamma_j}$ converge to the hyperbolic
structure of $M$.  In particular, there are holonomy characters
$\chi_j$ of the $M_{\gamma_j}$ which converge to $\chi^\df$.
Additionally, we can assume for each $j$ that the core of the solid
torus attached to $M$ to form $M_{\gamma_j}$ is a geodesic in
$M_{\gamma_j}$.  By Corollary~3.28 in \cite{Porti97}, $\chi^\df$ is a
smooth point of $X(M)$, and since the $\chi_j $ converge to
$\chi^\df$, infinitely many $\chi_j$ lie in $X_0$ (alternatively, for
our purposes, we could just change $X_0$ if necessary).

By Proposition~\ref{factor_theorem} the volume of a character depends
only on its image in $X(\bdry M)$, or more precisely in the smooth
projective model of $i^*(X_0)$.  This is key to the proof.

Consider one of the holonomy characters $\chi_j$ which comes from
Dehn-filling $M$ along the curve $\gamma_j$ in $\bdry M$.  From now
on, consider the map $i^*\maps X_0 \to i^*(X_0)$ as a map to a smooth
projective model of $i^*(X_0)$, though I will not change notation.  Let
$p_j = i^*(\chi_j)$.  Suppose $\chi$ is a point in $(i^*)^{-1} ( p_j)$
besides $\chi_j$.  Now if $\rho_j$ is a representation corresponding
to $\chi_j$ then $\rho_j(\gamma_j) = I$.  Let $\beta$ be a curve in
$\bdry M$ which forms a basis with $\gamma_j$ of $\pi_1(\bdry M)$.
Then $\beta$ is homotopic in $M_{\gamma_j}$ to the core of the solid
torus attached to $M$ to form $M_{\gamma_j}$, and so $\beta$ is homotopic
to a closed geodesic in $M_{\gamma_j}$.  So $\rho_j(\beta)$ is a
hyperbolic element of $\PSL{2}{\C}$.  In particular, $\tr(\rho_j(\beta))
\neq \pm 2$.  Now if $\rho$ is a representation whose character is
$\chi$ then $\tr(\rho(\gamma_j)) = \tr(\rho_j(\gamma_j)) = \pm 2$ and
$\tr( \rho(\beta) ) = \tr(\rho_j(\beta)) \neq \pm 2$.  Since
$\rho(\beta)$ and $\rho(\gamma_j)$ commute, we must have
$\rho(\gamma_j) = I$.  Hence $\rho$ is also a representation of
$\pi_1(M_{\gamma_j})$.

Since $\chi_j$ and $\chi$ map to the same point in $i^*(X_0)$, they
have the same volume.  Moreover by Lemma~\ref{factor_lemma}, for any
representation $\psi \maps \pi_1(M) \to \PSL{2}{\C}$ which factors
through $\pi_1(M_{\gamma_j})$, the volume of $\psi$ does not depend on
whether it is computed with respect to $M_{\gamma_j}$ or $M$.  Now by
Volume Rigidity for $M_{\gamma_j}$, the representation $\rho \maps
\pi_1(M_{\gamma_j}) \to \PSL{2}{\C}$ must be discrete and faithful.
Hence, as discussed in Section~\ref{df_reps}, $\rho$ is conjugate to
$\rho_j$ or $\bar{\rho}_j$.  I claim that for large $j$, $i^*(\rho_j)
\neq i^*(\bar{\rho}_j)$.  We know that the $\chi_j$ converge to
$\chi^\df$ and so the $\bar{\chi}_j$ converge to $\bar{\chi}^\df$.  By
Lemma \ref{distinct_images_of_d.f._reps} we know neighborhoods of
$\chi^{\df}$ and $\bar{\chi}^{\df}$ go to distinct branches of
$i^*(X(M))$.  So for large $j$, $i^*(\chi_j) \neq i^*(\bar{\chi}_j)$.
Hence for large $j$, $(i^*)^{-1} (p_j)$ consists only of $\chi_j$.  So
the map $i^* \maps X_0 \to i^*(X_0)$ has degree 1. \halmos
\end{proof}

It is easy to deduce the corresponding result for
$\SL{2}{\C}$ character varieties.

\begin{Corollary}\label{SL_degree_theorem}
Let $M$ be a finite-volume hyperbolic 3-manifold with one cusp.  Let
$\tilde{X}_0$ be a component of the $\SL{2}{\C}$ character variety of
$M$ which contains the character of a discrete faithful
representation.  The inclusion $i \maps \bdry M \injects M$ induces a
map $i^* \maps \tilde{X}_0 \to \tilde{X}(\bdry M)$.  This map has
degree onto its image at most $|H^1(M, \Z_2)|/2$ where $|\cdot |$
denotes number of elements.  In particular, if $H^1(M, \Z_2) = \Z_2$
then $i^*$ is a birational isomorphism onto its image.
\end{Corollary}

\begin{proof} 
If a representation $\rho$ into $\PSL{2}{\C}$ lifts to a
representation $\tilde{\rho}$ into $\SL{2}{\C}$ then there are
$|H^1(M, \Z_2)|$ distinct lifts which are constructed like this: If
$\epsilon \in H^1(M, \Z_2)$ is thought of as a homomorphism $\epsilon
\maps \pi_1 ( M )\to \Z_2 = \{I, -I \} \subset \SL{2}{\C}$ then we can
construct another lift of $\rho$ by $\phi(g) = \epsilon(g)
\tilde{\rho}(g)$.  By Poincar\'e duality and the long exact sequence
for the pair $(M, \bdry M)$, the image $H_1(\bdry M, \Z_2) \to H_1(M,
\Z_2)$ is one-dimensional.  So if $\chi$ is a $\PSL{2}{\C}$ character
which lifts to an $\SL{2}{\C}$ character, then the $|H^1(M, \Z_2)|$
distinct lifts map to precisely two points in $\tilde{X}(\bdry M)$,
unless the traces of $\pi_1(\bdry M)$ are all zero.  Since all the
traces are not zero generically on $\tilde{X}_0$, the last theorem
shows that the map $i^* \maps \tilde{X}_0 \to \tilde{X}(\bdry M)$ has
degree at most $|H^1(M, \Z_2)|/2$. \halmos
\end{proof}

\begin{Remark} 
In the first example of \cite{Dunfield2} the map $i^*$ has degree 4 on
$\tilde{X}_0$ and the first homology is $\Z \oplus \Z_4 \oplus \Z_2$.
For components of $X(M)$ which do not contain a discrete faithful
character, the map $i^*\maps X(M) \to X(\bdry M)$ may not have degree
1.  For instance, there are sometimes components of $X(M)$ which have
dimension greater than 1 whose image under $i^*$ is one-dimensional
(see Theorem~8.2 and 10.1 of \cite{CooperLong96}).  Even if we consider
an irreducible component of $X(M)$ of dimension 1 whose image under
$i^*$ is one-dimensional, $i^*$ may still fail to have degree 1.  This
happens with the exterior of the knot $7_4$.
\end{Remark}

\section{Cyclic Surgery Slopes and Boundary Slopes}\label{pf_of_full_theorem}

Let $K$ be a knot in a homotopy sphere $\Sigma$ with exterior $M$.
Fix a meridian $\mu$ and longitude $\lambda$ in $H_1(\bdry M, \Z)$.
The slope of $\gamma \in H_1(\bdry M, \Z)$ with respect to the basis
$(\mu, \lambda)$ will be denoted $r_\gamma$.  In this section I will
prove:

\begin{Theorem}\label{full_theorem}
Suppose $K$ is a small hyperbolic knot in a homotopy sphere.  Suppose
there is a $\beta \in H_1(\bdry M, \Z)$ with $\pi_1(M_\beta)$ cyclic
and $\beta \neq \pm \mu$.  Then there is an essential surface in the
exterior of $K$ whose boundary slope is non-integral and lies in
$(r_\beta - 1, r_\beta + 1)$.
\end{Theorem}

By the Cyclic Surgery Theorem \cite{CGLS}, $r_\beta$ is always an
integer.  Theorem~\ref{full_theorem} is true more generally for a
small hyperbolic knot in a manifold $\Sigma$ with cyclic fundamental
group whose exterior satisfies $H^1(M, \Z_2) = \Z_2$.  In this
setting, take $\lambda$ to be an arbitrary element of $H_1(\bdry M,
\Z)$ such that $(\mu, \lambda)$ is a basis (there is not always a
natural choice when $\pi_1(\Sigma)$ is non-trivial).  Note the
condition on $H^1(M, \Z_2)$ holds if $\pi_1(\Sigma)$ has odd order.
Theorem~\ref{full_theorem}, including this more general case, follows
easily from:

\begin{Theorem}\label{real_full_theorem}
Let $M$ be a one-cusped finite-volume hyperbolic 3-manifold with
$H^1(M, \Z_2) = \Z_2$.  Suppose $(\mu, \beta)$ is a basis for
$H_1(\bdry M, \Z)$ where $\pi_1(M_\beta)$ is cyclic.  For
$\gamma \in H_1(\bdry M, \Z)$, denote by $s_\gamma$ the slope with
respect to the basis $(\mu, \beta)$.  Then there is a boundary class
$\gamma$ such that $|s_\gamma| < 1$.
\end{Theorem}

Note that hypotheses of Theorem~\ref{real_full_theorem} do not
restrict $\pi_1(M_\mu)$, nor is it assumed that $M$ is small.  I will
now prove Theorem~\ref{full_theorem} assuming
Theorem~\ref{real_full_theorem}.

\begin{proof}[of Theorem~\ref{full_theorem}]
By the Cyclic Surgery Theorem, the meridian $\mu$ and the given class
$\beta$ form a basis for $H_1(\bdry M, \Z)$.  Let $\gamma$ be the
boundary class of $M$ given by Theorem~\ref{real_full_theorem}.  Note
$r_\gamma = s_\gamma + r_\beta$, and recall that $r_\beta$ is an
integer.  Since we are assuming $M$ is small, Theorem~2.0.3 of
\cite{CGLS} shows that $\beta$ is not a boundary class.  Therefore
$r_\gamma$ is a non-integral boundary slope in $(r_\beta -1, r_\beta +
1)$.
\end{proof}

I will now prove Theorem~\ref{real_full_theorem}:

\begin{proof} [of Theorem~\ref{real_full_theorem}]
I will use $\SL{2}{\C}$ character varieties here, so fix a
component $\tilde{X}_0$ of $\tilde{X}(M)$ which contains the
character of a discrete faithful representation.  I will break the proof
into three lemmas.  The first is:

\begin{Lemma}\label{ratio_constant}  
If $K$ is a counterexample to the theorem, then the function
$f_\mu/f_\beta$ is constant on $\tilde{X}_0$.
\end{Lemma}

So suppose $f_\mu/f_\beta$ is constant.  From this, it is possible to
determine $i^*(\tilde{X}_0)$ precisely.  To state the answer, I will
use the associated plane curve and the notation of
Section~\ref{plane_curve}.  Let $D_0$ be an irreducible component of
$t^{-1}(i^*(\tilde{X}_0))$.  Let $(m, b)$ be the eigenvalue
coordinates on $\Delta$ corresponding to the basis $(\mu, \beta)$.
The second lemma is:

\begin{Lemma}\label{deduce_D_0}
If $f_\mu/f_\beta$ is constant on $\tilde{X}_0$,  there is a constant
$C \neq \pm 1$ such that $D_0$ is exactly the set of zeros of the 
irreducible polynomial \begin{equation} \label{defining_p}
P(m, b) \equiv b m^2 - b - C b^2 m  + C m.
\end{equation}  
\end{Lemma}

The proof of the theorem is finished with:

\begin{Lemma}\label{finish_proof}
The polynomial $P$ of the last lemma can not define $D_0$ for a
hyperbolic manifold $M$ with $H^1(M, \Z_2) = \Z_2$.  
\end{Lemma}

I will now prove these lemmas.  The first follows immediately from the
following result, which is more general in that it makes no
assumption about $H^1(M, \Z_2)$:

\begin{Lemma} \label{special_slopes}
Let $M$ be a one-cusped finite-volume hyperbolic 3-manifold.  Suppose
$(\mu, \beta)$ is a basis for $H_1(\bdry M, \Z)$  where
$\pi_1(M_\beta)$ is cyclic.  Either there is a boundary class $\gamma$
such that $|s_\gamma| < 1$, or the function $f_\mu/f_\beta$ is
constant on $\tilde{X}_0$.
\end{Lemma}

\begin{proof}[Lemma~\ref{special_slopes}]
Set $g = f_\mu/ f_\beta$.  Let $Y$ be a smooth projective model of
$\tilde{X}_0$.  If $g$ is constant on $Y$, we are done.  Otherwise,
I will produce an essential surface with boundary class $\gamma$
where $|s_\gamma| < 1$.  As $s_\beta = 0$, if $\beta$ is a boundary
class take $\gamma = \beta$.  So assume $\beta$ is not a boundary
class. Let $y \in Y$ be a pole of $g$.  Suppose $f_\beta$ has a zero
at $y$.  For any rational function $h$ on $Y$, let $Z_y(h)$ denote the
order of zero at $y$, where $Z_y(h) = 0$ if $h$ does not have a zero
at $y$.  As $\beta$ is not a boundary class and $\pi_1(M_\beta)$ is
cyclic, Proposition~1.1.3 of \cite{CGLS} shows that the function
$f_\mu$ also has a zero at $y$ and $Z_y(f_\mu) \geq Z_y(f_\beta)$.
But then $g$ does not have a pole at $y$, a contradiction.  So
$f_\beta$ must not have a zero at $y$.  Thus $f_\mu$ must have a pole
at $y$, and so $y$ must be an ideal point where the associated
surfaces have non-empty boundary (see
Section~\ref{associated_surfaces}).  Let $\gamma$ be a boundary class
associated to $y$.  By the proof of Lemma~1.4.1 of \cite{CGLS}, we
have:
\begin{equation}\label{slope_at_ideal_point}
|s_\gamma| = \frac{\Pi_y(f_\beta)}{\Pi_y(f_\mu)},
\end{equation}
where $\Pi_y$ denotes the order of pole at $y$ or is $0$ if there is no
pole.  Since $\Pi_y(g) =  \Pi_y(f_\mu) - \Pi_y(f_\beta) > 0$, we have 
$|s_\gamma| < 1$.   \halmos
\end{proof}

Now I will prove Lemma~\ref{deduce_D_0}, which determines the equation
defining $D_0$, assuming that $f_\mu /f_\beta$ is constant.

\begin{proof}[Lemma \ref{deduce_D_0}]
Let $C'$ be the constant such that 
$f_\mu/f_\beta = C'$ on $D_0$.  With our coordinates $(m,b)$ on $\Delta$, 
we have 
\[
f_\mu = \left(m + \frac{1}{m} \right)^2 - 4
         = \left( m - \frac{1}{m} \right)^2 
\mtext{and}
f_\beta = \left(b - \frac{1}{b}\right)^2.
\]
So on $D_0$ the equation
\[
 \left(m - \frac{1}{m}\right)^2 = C' \left(b - \frac{1}{b}\right)^2
\]
holds.
 Since $D_0$ is irreducible, for some square root $C$ of $C'$ the 
 equation
\[
\left(m - \frac{1}{m}\right) = C \left(b - \frac{1}{b}\right)
\]
holds on $D_0$. Turning this into a polynomial condition, we have
\begin{equation} \label{new_defining_p}
P(m, b) \equiv b m^2 - b - C b^2 m  + C m = 0
\end{equation}
on $D_0$.  

Next, I will show $C \neq \pm 1$.   If $C = \pm 1$, $P = (b m + C)(m - C b)$.  
Then one of the eigenvalue functions $\xi_{\mu +\beta}$ or $\xi_{\mu - \beta}$ 
would be constant on $D_0$, which is impossible by Proposition~2 of
\cite{CullerShalen84}. 

If $C \neq \pm 1$, it is an elementary exercise to check that $P$ is irreducible.  
Thus $D_0$ must be {\em exactly} the zero set of $P$.
\halmos
\end{proof}

Finally, I will prove Lemma~\ref{finish_proof}, showing that the
polynomial $P$ of Lemma~\ref{deduce_D_0} can not be the defining
polynomial for $D_0$ because $H^1(M, \Z_2) = \Z_2$.

\begin{proof}[Lemma \ref{finish_proof}]
As before, let $Y$ be a smooth projective model of $\tilde{X}_0$.
Since $H^1(M, \Z_2) = \Z_2$, Corollary~\ref{SL_degree_theorem} shows
that the map $i^* \maps \tilde{X}_0 \to i^*(\tilde{X}_0)$ is a
birational isomorphism.  So we can also think of $Y$ as a smooth
projective model of $i^*(\tilde{X}_0)$.  Thus we have a rational map
$t \maps D_0 \to Y$ induced by $t \maps D_0 \to i^*(\tilde{X}_0)$.
Note $P(0,0) = 0$, and so $(0,0) \in D_0$.  This is a smooth point of
$D_0$ and the equation for the tangent line through $(0,0)$ is $C m -
b = 0$.  It follows that the eigenvalue functions $m = \xi_\mu$ and $b
= \xi_\beta$ restricted to $D_0$ have simple zeros at $(0,0)$.  So the
functions $\tr_\mu = m + 1/m$ and and $\tr_\beta = b + 1/b$ on $D_0$
have simple poles at $(0,0)$.  Let $y = t\left( (0,0) \right)$.  Then
$y \in Y$ is an ideal point of $\tilde{X}_0$ and both $\tr_\mu$ and
$\tr_\beta$ have simple poles at $y$.  Looking again at the line
tangent to $D_0$ though $(0,0)$, we see that $\mu - \beta$ is a
boundary class associated to $y$ (see
Section~\ref{associated_surfaces}).  Also, the eigenvalue of
$\mu - \beta$ at $(0,0)$ is $C$.  By the theorem in Section~5.7 of
\cite{CCGLS} the value of $\xi_{\mu - \beta}$ at $y$ is a root of unity
whose order divides the number of boundary components of any surface
associated to $y$.  So since $C \neq \pm 1$, any surface associated to
$y$ has at least three boundary components.  We now get our
contradiction by showing that there is a surface associated to $y$
with only two boundary components.  Let $T_y$ be the tree associated
to $y$.  By Proposition~\ref{associated_tree}, the translation length
of $\mu$ on $T_y$ is twice the order of pole of $\tr_\mu$ at $y$.
From the dual role of $Y$ as the smooth projective model of both
$i^*(\tilde{X}_0)$ and $\tilde{X}_0$, we have calculated that the order
of pole of $\tr_\mu$ at $y$ is 1, and so the translation length of
$\mu$ is 2.  But then Proposition~\ref{number_of_boundary_components}
shows there is a surface associated to $y$ with two boundary
components, a contradiction.  So $P$ can not be the equation defining
$D_0$. \halmos
\end{proof}

Since we have proved Lemmas \ref{ratio_constant}, \ref{deduce_D_0}, and
\ref{finish_proof}, we have proven the theorem. \halmos
\end{proof}

\begin{Remark} 
Note that the proof used $H^1(M, \Z_2) = \Z_2$ in a fundamental way.
If the degree of $\tilde{X}_0 \to i^*(\tilde{X}_0)$ were $d$,
we could only have concluded that there was a surface associated to
$y$ having as few as $2d$ boundary components.  The first example in
\cite{Dunfield2} has an ideal point whose image in
$i^*(\tilde{X}_0)$ looks exactly like in the proof, in the sense
that $\tr_\mu$ has a simple pole at the corresponding point of
$D_0$.  But there,  $H_1(M, \Z)$ is $\Z \oplus \Z_4 \oplus \Z_2$ and the
associated surface has four boundary components.  This illustrates
the way that Corollary~\ref{SL_degree_theorem} shows there is a 
surprising connection between the size of $H^1(M, \Z_2)$ and the
character variety of $M$.
\end{Remark}

\section{The norm and the diameter of the set of boundary slopes
}\label{norm}

Lemma~\ref{special_slopes} has an interpretation in terms of the norm
on $H_1(\bdry M, \R)$ defined in \cite{CGLS}.  Let $V$ denote
$H_1(\bdry M, \R)$ and $L \subset V$ the lattice $H_1(\bdry M, \Z)$.
This norm on $V$ has the property that for $\gamma \in L$, the number
$||\gamma||$ is the degree of $f_\gamma$ on $\tilde{X}_0$. Let $r$ be
the minimum of $||\gamma||$ over all non-zero $\gamma \in L$.  Let $B$
denote the closed ball of radius $r$ about $0$ in $V$.  The ball $B$ is
a finite-sided convex polygon which is invariant under $v \mapsto -v$.
Boundary classes associated to ideal points of $\tilde{X}_0$ correspond
bijectively with vertices of $B$, where a boundary class $\gamma$
corresponds to a vertex $v$ where $v = a \gamma$ for a positive
rational number $a$ (this is not quite explicit in \cite{CGLS}, but see
Lemma~6.1 of \cite{BoyerZhang96}). An element $\gamma
\in L$ such that $\pi_1(M_\gamma)$ is cyclic and which is not a boundary
class has minimal norm, that is, $||\gamma|| = r$ and $\gamma \in B$.
In this language, a slight modification of the proof of
Lemma~\ref{special_slopes} gives:

\begin{Lemma}\label{modified_lemma}
Suppose $(\mu, \beta)$ is a basis of $L$ such that
$\pi_1(M_\mu)$  is cyclic, $\mu$ is not a boundary slope, and $\beta$
has minimal norm. Then either there is a boundary class $\gamma$,
associated to an ideal point of $\tilde{X}_0$, which satisfies $|s_\gamma|< 1$
or the function $f_\mu/f_\beta$ is constant on $\tilde{X}_0$. 
\end{Lemma}

\begin{proof}  
The needed modification to the proof of Lemma~\ref{special_slopes} is
at the step that shows that if $f_\beta$ has a zero at $y$ then
$f_\mu$  also has a zero at $y$ and $Z_y(f_\mu) \geq
Z_y(f_\beta)$.  Here, by Proposition~1.1.3 of \cite{CGLS}, we know if 
$f_\mu$ has a zero at $y$ then $Z_y(f_\mu) \leq Z_y(f_\beta)$.  
The total number of zeros of $f_\mu$
is equal to the total number of zeros of $f_\beta$ since $||\mu|| =
||\beta||$.  Thus the sets of zeros
(including multiplicity) of $f_\mu$ and $f_\beta$ are the same, and we
can apply the rest of the proof of Lemma~\ref{special_slopes} unchanged.
 \halmos
\end{proof}

I will now give an application to the following question.  
Consider a knot in a homotopy sphere with irreducible
exterior $M$. The set of boundary slopes is finite \cite{Hatcher82},
and so has a well-defined diameter $d$ as a subset of $\Q \cup \{
\infty \}$ (I will use the convention that $d = \infty$ if $\infty$ is a
boundary slope).  Hatcher and Thurston \cite{HatcherThurston} asked
in the case of knots in $S^3$ whether $d$ is always greater than $2$. In
\cite{CullerShalenDiameter}, Culler and Shalen showed that for any
knot, the diameter $d \geq 2$.  Another application of the proof of
Theorem~\ref{full_theorem} is to show:

\begin{Theorem}
If $d = 2$ for a hyperbolic knot in a homotopy sphere, the greatest
and least boundary slopes are not integers.
\end{Theorem}

\begin{proof}
Consider such a knot with $d = 2$.  In \cite[Proof of Main
Theorem]{CullerShalenDiameter} it is shown that for a suitable choice
of meridian-longitude basis $(\mu, \lambda)$ for $L$, $B$ is as in
Fig.~\ref{fund_poly}.  Since this is not quite explicit there, let me
elaborate.  In the notation of \cite{CullerShalenDiameter}, let $v_1$
and $v_2$ be the vertices of $B$ which are the ends of the edge
containing $\mu$.  Then $\mu = t v_1 + (1 - t) v_2$ for some $t \in
[0,1]$, and they show that the diameter of the set of boundary slopes is
bounded below by $1/(2t(t-1)) \geq 2$.  If the diameter of the set of
boundary slopes is 2, we must have $t = 1/2$.  In this case, the area
of the parallelogram with vertices $\{v1, v2, -v1, -v2\}$ is 4 and so
$B$ is equal to this parallelogram.  Since $t = 1/2$, $\mu$ lies in
the middle of an edge of $B$.  Combined with the fact that the area of
$B$ is $4$, the sides of $B$ must be segments of vertical lines
$\lambda = \pm 1$.  Since $0$ is always a boundary slope, after
possibly changing the signs of $\mu$ and $\lambda$, $B$ must be as
in Fig.~\ref{fund_poly}.

\begin{figure}[hbt]
  \centerline{\psfig{figure=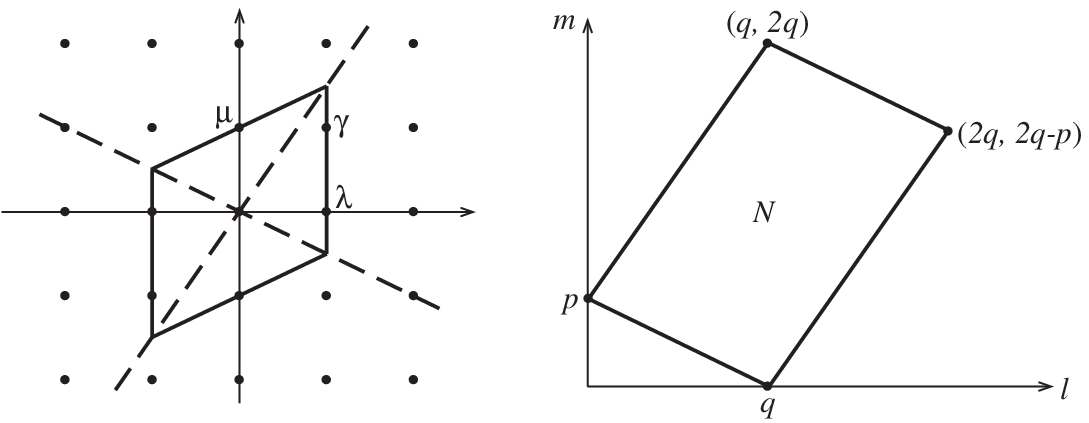}}
  \firstcaption{The fundamental polygon.  The slopes of the dotted
    lines are $2 - p/q$ and $-p/q$.          \label{fund_poly}}   
    \secondcaption{The Newton polygon of $P$ \label{Newton_polygon}}
\end{figure}
From the correspondence between boundary classes associated to ideal
points of $\tilde{X}_0$ and vertices of $B$, there are only two
boundary slopes associated to ideal points of $\tilde{X}_0$, namely,
$-p/q$ and $2 - p/q$ where $0 \leq p \leq q$ and $\gcd(p,q) = 1$.

Showing that the greatest and least slopes are not integers is
equivalent to showing $q > 1$.  Suppose $q = 1$.  There are two cases:
$p = 0$ and $p = 1$.  If $p = 1$, the two boundary slopes are $1$ and
$-1$ with respect to the basis $(\mu, \lambda)$.  Since $\lambda$ is in
$B$ it has minimal norm, and as $\mu$ is not a boundary class (since
$d$ is finite), Lemma \ref{modified_lemma} shows that
$f_\mu / f_\lambda$ must be constant.  We can now apply
Lemmas~\ref{deduce_D_0} and \ref{finish_proof} to get a contradiction.
If $p = 0$, we apply the same argument with $\lambda$ replaced by
$\gamma = \mu + \beta$.  \halmos
\end{proof}

It is possible to prove more:

\begin{Theorem}\label{diameter}
Suppose the diameter of the set of boundary slopes of a hyperbolic
knot in a homotopy sphere is 2.  Let $0 \leq p \leq q$ be such that
the greatest and least slopes are $2-p/q$ and $-p/q$ respectively.
Then $p$ is even, $q$ odd, and $q > 1$.
\end{Theorem}

The key to this  version is a lemma which is a consequence of
Theorem~\ref{degree_theorem}.  Consider the commutative diagram:
 \begin{equation}
\xymatrix{  {\tilde{X}_0} \ar[r]^{i^*} \ar[d]_{\pi_M} 
 & **[r] i^*(\tilde{X}_0) \subset \tilde{X}(\bdry M) \ar[d]^{ \pi_{\bdry M}} \\
 X_0 \ar[r]^{i^*} & **[r] i^*(X_0) \subset X(\bdry M)
}
\label{diagram}
\end{equation}
The map $\pi_{\bdry M}$ from $\tilde{X}(\bdry M)$ to $X(\bdry M)$ has
degree 4, so one might expect that the restriction of 
$\pi_{\bdry M}$ to $i^*(\tilde{X}_0)$ could also have degree 4.  But I will
show the degree of this restriction can be no more than 2.   An
involution $f$ of  $i^*(\tilde{X}_0)$ which has the property that
$\pi_{\bdry M} \circ f = \pi_{\bdry  M}$ will be called a {\em symmetry}
of  $i^*(\tilde{X}_0)$.   If the restriction of $\pi_{\bdry M}$ 
to $i^*(\tilde{X_0})$ has degree 4, the group of symmetries would 
have order 4.
One possible symmetry, $\tau$, is the restriction of the involution of
$\tilde{X}(\bdry M)$ whose action in coordinates is 
$( \tr_\mu, \tr_\lambda, \tr_{\mu +
\lambda} ) \mapsto (-\tr_\mu, \tr_\lambda, - \tr_{\mu + \lambda})$.  I will show:

\begin{Lemma}\label{symmetry}
   For the exterior of a hyperbolic knot in a homotopy sphere, the only
   possible symmetry of $i^*(X_0)$ is $\tau$.
 \end{Lemma} 

\begin{proof}
The discussion in the proof of Corollary~\ref{SL_degree_theorem} shows
that $\pi_M$ has degree 1 or 2, depending on whether
$(\pi_M)^{-1}(X_0)$ is irreducible.  By Theorem~\ref{degree_theorem}
and Corollary~\ref{SL_degree_theorem}, the horizontal maps in
Eqn.~(\ref{diagram}) have degree 1.  So the degrees of $\pi_M$ and
$\pi_{\bdry M}$ are equal.  If these degrees are 1, there are no
symmetries.  If the degrees are 2, $i^*(\tilde{X}_0)$ has one
symmetry, which I claim is $\tau$.  In this case,
there is an involution $\tau'$ of $\tilde{X}_0$ for which $\pi_M \circ
\tau' = \pi_M$, namely, multiplication of characters by the unique
non-trivial homomorphism $\pi_1(M) \to \Z_2 = \{I, -I \} \subset
\SL{2}{\C}$.  Then $\tau'$ induces the symmetry $\tau$ of
$i^*(\tilde{X})$. \halmos
\end{proof}  

\begin{Remark}  For the exterior of a knot in a homotopy sphere, it
is conceivable that $i^*(X_0)$ might not have the symmetry $\tau$.
Suppose $(\pi_M)^{-1}(X_0)$ consists of two irreducible components.
While $i^*\left( (\pi_M)^{-1}(X_0) \right)$ would be invariant under
$\tau$, the curve $i^*(\tilde{X}_0)$ might not be.  The first example
of \cite{Dunfield2} is a one-cusped hyperbolic manifold where
$(\pi_M)^{-1}(X_0)$ splits into multiple components.  Applying $i^*$ to
one component yields a curve which is not invariant under the analogue
of $\tau$.  However, this manifold is not exterior of a knot in a
homotopy sphere.
\end{Remark}
 
From Lemma~\ref{symmetry} we can prove Theorem~\ref{diameter}:
 
\begin{proof}[Theorem~\ref{diameter}]

Let $p$ and $q$ be as in the statement.  Theorem~\ref{diameter} shows
that $q > 1$, and I will assume this throughout.  Let $\gamma = \mu + \lambda$.
Consider the function
\[
g = \frac{f_\lambda^p f_\gamma^{q - p}}{f_\mu^{q}} \mbox{\ on\ }
\tilde{X}_0,
\]
noting that $q - p > 0$. The proof follows the same basic plan as that of
Theorem~\ref{full_theorem}.  The analogue of Lemma~\ref{ratio_constant} is:

\begin{Lemma}\label{new_g_const}
If $d = 2$, then $g$ is constant on $\tilde{X}_0$. 
\end{Lemma}

Let $(m, l)$ be the eigenvalue coordinates on $\Delta$ corresponding
to $(\mu, \lambda)$.  The analogue of Lemma~\ref{deduce_D_0} is:

\begin{Lemma}\label{new_deduce_D_0} If $g$ is constant
there is a $C \in \C$ so that $D_0$ is exactly the zeros of
\begin{equation}
P(m,l) \equiv m^p \left(l^2 - 1\right)^{p} \left(l^2 m^2 - 1\right)^{q
  - p} - C l^q \left(m^2 -
1 \right)^q. \label{new_P} \end{equation}
\end{Lemma}
 
The final step, which is the one that differs most from the proof of
Theorem~\ref{full_theorem},  is to use Lemma~\ref{symmetry} to show the
following analogue of Lemma~\ref{finish_proof}:

\begin{Lemma}\label{new_proof}
If the polynomial $P$ of the last lemma defines $D_0$ for the exterior
of a knot in a homotopy sphere, then $p$ is even and $q$ is odd.
\end{Lemma}

I will begin with:

\begin{proof}[Lemma \ref{new_g_const}]
This is a refinement of the proofs of Lemmas~\ref{special_slopes} and
\ref{modified_lemma}.  Suppose $g$ is not constant.  Let $y \in Y$ be
a pole of $g$, where $Y$ is a smooth projective model of
$\tilde{X}_0$.  We know that $\pi_1(M_\mu)$ is a cyclic and that $\mu$
is not a boundary class.  From Fig.~\ref{fund_poly} we know $\mu$,
$\lambda$, and $\gamma$ all have minimal norm.  Therefore, as in the
proof of Lemma~\ref{modified_lemma}, the sets of zeros (including
multiplicity) of $f_\mu$, $f_\lambda$, and $f_\gamma$ are all the
same.  But then $g$ would not have a pole at $y$, a contradiction.
Thus at least one of $f_\lambda$ or $f_\gamma$ has a pole at $y$.  So
$y$ is an ideal point.  

By Lemma 1.4.1 of \cite{CGLS} there is a linear functional $l$ on
$H_1(\bdry M, \Z)$ such that the order of pole of $f_\alpha$, $\Pi_y
f_\alpha$, is $| l(\alpha) |$ for all $\alpha$ in $H_1(\bdry M, \Z)$.
The slope associated to $y$ is either $-p/q$ or $2 - p/q$.  In the
first case, $\Pi_y f_{-p\mu + q\lambda} = 0$ and so there is a $d \in
\Z$ such that $l( a \mu + b \lambda) = d ( qa + pb)$.  Then $\Pi_y
f_\mu = |d| q, \Pi_y f_\lambda = |d| p$, and $\Pi_y f_\gamma = |d| (p
+ q)$.  So 
\[ \Pi_y g = -q \Pi_y f_\mu + p \Pi_y f_\lambda + (q - p)\Pi_y f_\gamma =
0.\]  
If the slope associated to $y$ is $2 - p/q = (2q - p)/q$ then $l( a
\mu + b \lambda) = d ( qa - (2q - p)b)$.  Then $\Pi_y f_\mu = |d| q,
\Pi_y f_\lambda = |d|(2q - p)$, and $\Pi_y f_\gamma = |d| (q- p)$;
hence $\Pi_y g = 0$.

So $g$ has no poles and must be constant.  \halmos
\end{proof}

Now we deduce the equation defining $D_0$.  

\begin{proof}[Lemma \ref{new_deduce_D_0}]

Just as in the proof of Lemma~\ref{deduce_D_0}, there is a
constant $C$ so that the following equation holds on $D_0$:
\begin{equation}\label{eqn_def_new_D_0}
\left(l - \frac{1}{l}\right)^p \left(ml- \frac{1}{ml}\right)^{q-p}
= C \left(m - \frac{1}{m}\right) ^q,
\end{equation}
This is equivalent to the polynomial $P$ given in the statement of
Lemma~\ref{new_deduce_D_0} being zero.

Now I will show that $D_0$ is exactly the zeros of $P$.  
We need to take the point of view of \cite{CCGLS} in which information
about ideal points of $\tilde{X}_0$ is deduced from the Newton polygon of the
equation defining $D_0$.  The {\em Newton polygon} of a polynomial $Q$
in variables $x$ and $y$ is the convex hull in $\R^2$ of:
\[
 \{ (i, j) \in \Z^2 \ |  \ \mbox{the coefficient of $x^i y^j$ in $Q$ is
 nonzero}    \}.
\]
The Newton polygon, $N$, of $P$ is shown in Fig.~\ref{Newton_polygon}.

I claim that $\{ P = 0 \}$ is exactly $D_0$.  Suppose $Q$ were a
proper factor of $P$ which defines $D_0$.  Let $N'$ be the Newton
polygon of $Q$.  In \cite{CCGLS} it is shown that the slopes of the
sides of $N'$ are precisely the boundary slopes of surfaces
associated to ideal points of $\tilde{X}_0$.  Thus the slopes of the
sides of $N'$ are $-p/q$ and $2 - p/q$.  So we know that $N'$ is a
parallelogram whose sides are parallel to those of $N$.  The idea is
that $N$ is the smallest such parallelogram and therefore $N = N'$ and
$P = Q$.  Let $\diam_l (N')$ denote the diameter of the projection of
$N'$ onto the axis corresponding to the exponent of $l$.  If we think
of $Q$ as the a polynomial in a single variable $l$ over $\C[m]$, then
the degree of this polynomial, $\deg_l(Q)$, is equal to $\diam_l(N')$.
Since $Q$ is a factor of $P$, we know $\deg_l Q \leq \deg_l P$, and so
$\diam_l(N') \leq \diam_l(N) = 2q$.  Projecting onto the other axis,
we also have $\diam_m(N') \leq \diam_m(N) = 2q$.  Moreover, if both
diameters of $N'$ and $N$ agree, we have $P = Q$.  Now consider a side
$S$ of $N'$ with slope $-p/q$.  Because the endpoints of $S$ are in
$\Z^2$, we must have $\diam_l(S) \geq q$ and $\diam_m(S) \geq p$.
Similarly, a side $T$ of $N'$ with slope $2 - p/q$ must have
$\diam_l(T) \geq q$ and $\diam_m(T) \geq 2q - p$.  Thus $\diam_l (N')
\geq 2q = \diam_l(N)$ and $\diam_m(N') \geq 2q = \diam_m(N)$.  So $P =
Q$ as desired.  So $D_0$ must be exactly the zero set of $P$.  
\halmos 
\end{proof}

Now I will show that if $P$ defines $D_0$ for the exterior of a knot
in a homotopy sphere, then $p$ is even and $q$ is odd.  By Lemma
\ref{symmetry} the only allowed symmetry is the one whose action on
$D_M$ sends $(m, l)$ to $(-m,l)$.  If $q$ is even, $p$ and $q - p$ are
odd.  Note in this case that Eqn.~(\ref{eqn_def_new_D_0}) and hence
$P$ are invariant under $(m, l) \mapsto (m, -l)$.  But then there is a
symmetry of $i^*(\tilde{X}_0)$ other than $\tau$, which is impossible.
Therefore $q$ is odd, and we already know $q \neq 1.$ If $p$ and $q$
are both odd, then Eqn.~(\ref{eqn_def_new_D_0}) is invariant under
$(m, l) \mapsto (-m, -l)$, so this case is ruled out as well.  So $p$
is even and $q$ is odd. \halmos
\end{proof}

\begin{Remarks}

Let $Q$ denote the polynomial in Eqn.~(\ref{new_P}) with $p = 1$, $q =
2$, and $C = 1$.  The polynomial $Q$ actually occurs as the defining
equation for $D_0$ for $N$, the sister of the exterior of the
figure-8 knot. In this case, $N_\mu$ has fundamental group $\Z_{10}$.
The results of \cite{CullerShalenDiameter}, more generally, say that
$d \geq 2$ for any knot in a manifold with cyclic fundamental group
whose exterior is irreducible and not cabled.  It turns out that for
$N$, the diameter of the set of all boundary slopes is exactly $2$,
and this shows that the estimate of Culler and Shalen is sharp in this
more general context (see Example 1.4 of \cite{CullerShalenDiameter}).
It was the observation that equation defining $D_N$ was $\{ Q = 0 \}$
that led me to discover this example.

If one cares about the diameter of the set of {\em strict} boundary
slopes, then the situation stays the same except that $\lambda$ is
replaced by some random class $\nu$ with integer slope.  In this case,
you can not rule out $p$ and $q$ both being odd since $\nu$ may generate
$H_1(M, \Z_2)$.  This is also why you can not use the argument about
symmetries to prove Theorem~\ref{full_theorem}.

\end{Remarks}

\section{Proof of volume rigidity} \label{pf_of_volume_theorem}

This section provides a proof of:

\begin{Theorem}[Gromov-Thurston-Goldman]\label{volume_rigidity}
Suppose $M$ is a compact hyperbolic 3-manifold.  If $\rho_1 \maps
\pi_1 (M) \to \PSL{2}{ \C}$ is a representation with $\vol(\rho_1) =
\vol(M)$, then $\rho_1$ is discrete and faithful.
\end{Theorem}

\begin{proof}
The proof is essentially the same as that of Thurston's strict version
of Mostow's Theorem \cite[Theorem 6.4]{ThurstonLectureNotes}.  The
only modification is that since $\H^3/\rho_1$ may be nasty, rather
than a compact manifold, it is necessary to do some things
equivariantly.  I will follow \cite{ThurstonLectureNotes}, with some
details coming from Toledo's paper \cite{Toledo} using the same
technique.  The same proof works in higher dimensions with the aid of
\cite{HaagerupMunkholm}, but I stick to the 3-dimensional case for
simplicity.  I assume some familiarity with Gromov's proof of Mostow's
Theorem (for a nice account, see \cite{Munkholm} or the very through
\cite[Chapter 11]{Ratcliffe}).

Let $\rho_0$ be a discrete faithful representation for $M$.  Pick a
smooth equivariant map $f$ from $\H^3$ acted on by $\rho_0$ to $\H^3$
acted on by $\rho_1$.  If the integral $\int_M f^*(\Vol_\H^3)$ is negative,
choose an orientation reversing isometry $r$ of $\H^3$ and replace
$\rho_1$ by $r \circ \rho_1 \circ r^{-1}$ and $f$ by $r \circ f$ so that
the integral is positive.  Then the hypothesis on $\vol(\rho_1)$ gives:
\begin{equation}
 \int_M \Vol_{\H^3} = \vol(M) = \vol(\rho_1) = \int_M f^*(
 \Vol_{\H^3}). \label{condition}
\end{equation}

By a {\it tetrahedron}, I will mean a simplex in $\H^3$ with totally
geodesic faces.  I will divide the proof into the following three claims:

\begin{Claim}\label{f_extends}
The map $f$ extends to an equivariant measurable
map $\bar{f}$ from $S_\infty^2 = \bdry \H^3$ to itself.
\end{Claim}

\begin{Claim}\label{pres_reg_tet}
The map $\bar{f}$ takes vertices of almost all regular ideal
tetrahedra to vertices of regular ideal tetrahedra.
\end{Claim}

\begin{Claim}\label{done}
Because of Claim~\ref{pres_reg_tet}, $\bar{f}$ is essentially a
M\"obius transformation.  That is, there is a M\"obius transformation
$F$ so that $F = \bar{f}$ almost everywhere.
\end{Claim}

The proof is then completed by noting that since $F$ is equivariant,
it conjugates the action of $\rho_0$ on $S^2_\infty$ to the action of
$\rho_1$ on $S^2_\infty$.  As $\rho_0$ and $\rho_1$ are conjugate,
$\rho_1$ is discrete and faithful.

Let me start in on the proof of Claim~\ref{f_extends}.  The key is
that the regular ideal tetrahedron is the unique tetrahedron of
maximal volume \cite[Theorem 10.4.7]{Ratcliffe}.  The idea behind
Claim~\ref{f_extends} is that $f$ must take the vertices of a
non-ideal tetrahedron of near-maximal volume to points which span a
tetrahedron of near-maximal volume or else $f$ would be volume
shrinking and Eqn.~(\ref{condition}) would be violated.

The proof will use measure homology, which I will briefly describe
(for details see Section 4 of \cite{Munkholm} or Section 11.5 of
\cite{Ratcliffe}).  Let $C^1( \Delta_k, M)$ be the space of $C^1$ maps
of the standard $k$-simplex $\Delta_k$ into $M$, supplied with the
$C^1$ topology.  Let $\sC_k(M, \R)$ be the set of real-valued Borel
measures of bounded total variation on $C^1(\Delta_k, M)$.  There is a
natural boundary operator from $\sC_k(M, \R)$ to $\sC_{k-1}(M, \R)$,
and the homology of the resulting complex is called the {\it measure
  homology\/} of $M$.  There is an inclusion of the usual $C^1$
singular chain complex $C_*(M, \R)$ into $\sC_*(M, \R)$ which
sends a singular simplex (an element of $C^1(\Delta_k, M)$) to the
associated Dirac measure.  This inclusion induces an isomorphism on
homology.  The pairing between $C_{k}(M, \R)$ and smooth differential
$k$-forms extends to $\sC_k(M, \R)$.

Let $G = \PSL{2}{\C}$ and $\Gamma \subset G$ be $\pi_1(M)$ acting via
$\rho_0$.  The unit tangent bundle to $M$ is $X = \Gamma \backslash
G$, and let $\mu$ be the Haar measure on $G$ such that $\mu(X) =
\vol(M)$.  Let $D$ be a polyhedron minus some faces which is a
fundamental domain for $M$ in $\H^3$ (each orbit of $\H^3$ under
$\Gamma$ has exactly one representative in $D$).

Let $\sigma \subset \H^3$ be a non-ideal tetrahedron.  Let $\smear
\sigma$ be the measure cycle in $\sC_3(\H^3, \R)$ consisting of all
translates of $\sigma$ with first vertex in $D$, uniformly weighted by
$1/\vol(\sigma)$; that is, if $\eta$ is a differential 3-form:
\[
\pair{\eta, \smear \sigma} = \int_{g \in X} \left(
\frac{1}{\vol(\sigma)} \int_{g \cdot \sigma} \eta \right) \ d\mu.
\]
Here, I abuse notation and denote $D \cross \SO(3) \subset G$, a
fundamental domain for the unit tangent bundle $X$, by $X$.  Now if
$\sigma^-$ denotes $\sigma$ with the opposite orientation, let
$z_\sigma = 1/2( \smear \sigma - \smear \sigma^-)$.  If $\pi:\H^3 \to
M$ is the covering map, we have $\pi_*(z_\sigma)$ is closed.  Since
measure homology is the same as standard homology, $z_\sigma$
represents some multiple of the fundamental class.  As $\pair{
  \Vol_{\H^3}, \smear \sigma} = \vol(M)$, the cycle $z_\sigma$
represents the fundamental class.  Now $\vol(\rho_1) = \pair{f^*
  \Vol_{\H^3}, z_\sigma} = \pair{f^* \Vol_{\H^3}, \smear \sigma}$, and
so
\begin{equation}\label{replace}
\vol(\rho_1) = \int_X \left( \frac{1}{\vol(\sigma)} \int_{g \cdot
  \sigma} f^* \Vol_{\H^3} \right) \ d\mu.  
\end{equation}
There is a chain map $\Str$ from $\sC_*(\H^3, \R)$ to itself which
replaces a singular simplex by a geodesic simplex with the same
vertices (see Section 4 of \cite{Munkholm} or Section 11.5 of
\cite{Ratcliffe} for a definition of $\Str$ and its properties).  I
claim that we can replace $ \int_{g \cdot \sigma} f^* \Vol_{\H^3}$ in
Eqn.~(\ref{replace}) by $\int_{\Str(f(g \cdot \sigma))} \Vol_{\H^3} =
\vol(\Str( f( g \cdot \sigma)))$.  The idea is this: Consider the
analogous question for singular simplicial homology, i.e.~let $z$ be a
lift to $C_3(\H^3, \Z)$ of a cycle in $C_3(M, \Z)$.  The chain $z$ is
now a finite linear combination of singular simplices.  While $\bdry
z$ is not zero, you can choose elements of $\Gamma$ that pair up the
simplices of $\bdry z$ that reflect the fact that $\pi(\bdry z) = 0$
in $C_2 (M, \Z)$.  The difference between $f_*(z)$ and $\Str f_*(z)$
depends on $f_*(\bdry z)$.  The way elements of $\Gamma$ pair up the
simplices of $\bdry z$ shows that we can replace $ \int_{g \cdot
  \sigma} f^* \Vol_{\H^3}$ by $\int_{\Str(f(g \cdot \sigma))}
\Vol_{\H^3}$ in (\ref{replace}).

Formally, there is a chain homotopy, $H$, from $\Str$ to the identity
which is invariant under isometries.  Now
\[
\pair{f^*\Vol_{\H^3}, z_\sigma} = \pair{\Vol_{\H^3}, f_*(z_\sigma)}
= \pair{ \Vol_{\H^3}, \Str(f_*(z_\sigma)) - H(f_*(\bdry
  z_\sigma))}
\]
and so to show the claim it is enough to show $\pair{\Vol_{\H^3},
H(f_*(\bdry
  z_\sigma))} = 0$.  Equivalently,  define a cochain $c$ by
$c(\tau) = \pair{\Vol_{\H^3}, H(f_* (\tau))}$; we want $c(\bdry
z_\sigma) = 0$.  Since $f$ is equivariant and $H$ commutes with
isometries, $c$ descends to a cochain $c_M$ on $M$.  As $\pi_*(\bdry
z_\sigma) = 0$, $c_M(\pi_*(\bdry
z_\sigma)) = c(\bdry z_\sigma) = 0$, as desired.  This proves the
claim.  Hence we have from Eqn.~(\ref{replace}):
\begin{equation}
\vol(\sigma) \vol(M) = \int_X \vol(\Str( f( g \cdot \sigma ))) \
d\mu. \label{key}
\end{equation}
where the volume of $\Str( f( g \cdot \sigma ))$ is signed volume. This
is the formula which will guarantee that $f$ do not shrink large
tetrahedra very much.

Fix a geodesic ray $r$ with endpoint $b \in D$.  Let $\sigma_i$ be a
regular tetrahedron all of whose sides have length $i$ with one vertex
$b$ and an edge lying on $r$. The idea is to use the expanding
sequences $\{g \cdot \sigma_i \}$ for $g \in X$ to approximate what
would happen to an ideal tetrahedron.  Let $v_3$ denote the volume of
a regular ideal tetrahedron, which is the unique tetrahedron of
maximal volume.  The next lemma is a quantitative version of the
statement ``$f$ does not shrink volume.''

\begin{Lemma}\label{estimation} Let $\epsilon_i = v_3 -
  \vol(\sigma_i)$, and let $Y = \{g \in X : \vol(\Str( f(g \cdot
  \sigma_i)) < \vol(\sigma_i) - i^2 \epsilon_i \}$.  Then $\mu(Y) <
  \mu(X)/i^2$.
\end{Lemma}

\begin{proof} This is Lemma 2.3 of \cite{Toledo}.  Since 
  $\vol(\Str( f( g \cdot \sigma)) ) < v_3 = \vol(\sigma_i) +
  \epsilon_i$, Eqn.~(\ref{key}) gives us 
  \begin{eqnarray*}
  \vol(\sigma_i) \vol(M)
 &=& \int_X \vol(\Str( f (g \cdot \sigma_i ))) \ d\mu \\
 &=& \int_Y \vol(\Str( f (g \cdot \sigma_i ))) \ d\mu + \int_{X \setminus
 Y}   \vol(\Str( f ( g \cdot \sigma_i ))) \ d\mu \\
 &<& (\vol(\sigma_i) - i^2
  \epsilon_i)\mu(Y) + v_3  \mu(X \setminus Y)
\end{eqnarray*}
Since 
\[
(\vol(\sigma_i) - i^2 \epsilon_i)\mu(Y) + v_3 \mu(X \setminus Y) =
\vol(\sigma_i) \vol(M) + \epsilon_i (\mu(X \setminus Y) - i^2 \mu(Y))
\]
we have
\[
\vol(\sigma_i) \vol(M) < \vol(\sigma_i) \vol(M) + \epsilon_i (\mu(X
\setminus Y) - i^2 \mu(Y)).
\]
Therefore $0 < \mu(X \setminus Y) - i^2 \mu(Y)$, and so $\mu(Y) <
\mu(X)/i^2$. \halmos
\end{proof}

The next lemma shows that for large $i$, $\sigma_i$ is a very good 
approximation of a regular ideal tetrahedron. 

\begin{Lemma}\label{hyperbolic} Let $\sigma_i$ be a regular tetrahedron
in $\H^3$ all of whose sides have length $i$.  Then for $i$ large,
$\epsilon_i = v_3 - \vol(\sigma_i)$ decreases exponentially with $i$.
\end{Lemma}

\begin{proof}  This is Lemma 6.4.1 of \cite{ThurstonLectureNotes}.
Let $\delta_\infty$ be a fixed regular ideal tetrahedron.  Let $p$ be
the center of mass of $\delta_\infty$ and consider the four rays
starting at $p$ and ending at the vertices of $\delta_\infty$.  The
tetrahedron whose vertices are the points on these rays a distance $t$
from $p$ is regular, and will be denoted $\delta_i$ where $i \in
\R^+$ is the length of any of its edges (thus $\sigma_i$ and
$\delta_i$ are isometric).
\begin{figure}[ht]
\centerline{\psfig{figure=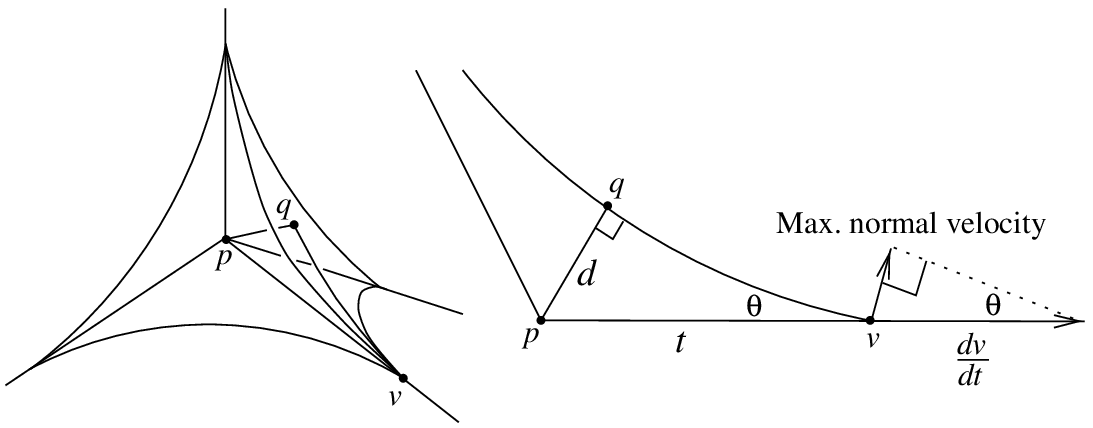}}
\firstcaption{The tetrahedron $\delta_i$\label{normal_velocity1}}
\secondcaption{Finding the maximum normal
  velocity\label{normal_velocity2}}
\end{figure}
Any tetrahedron has a natural straight (or barycentric)
parameterization coming from the affine structure of the hyperboloid
model of $\H^3$ (see Section 11.4 of \cite{Ratcliffe}).
Parameterizing $\delta_i$ in this way, we get a one parameter family
of diffeomorphisms $\varphi_i$ from the standard 3-simplex $\Delta_3$
to $\delta_i$.  A point $y = \varphi_i(x)$ on a face of $\delta_i$ has
a velocity which is defined as the derivative of $\varphi_i(x)$ with
respect to $i$.  The normal velocity of $y$ is the component of
velocity normal to the face of $\delta_i$ containing $y$.

The derivative $\frac{d \mbox{\small vol}(\delta_i)}{dt}$ is bounded
by the area of $\bdry \delta_i$ times the maximum normal velocity of
$\bdry \delta_i$.  Let $q$ be the center of mass of one of the faces
of $\delta_i$.  From Figs.~\ref{normal_velocity1} and
\ref{normal_velocity2}, we see that the maximum normal velocity of
$\bdry \delta_i$ is $\sin \theta$.

By the law of sines, $\sin \theta = \frac{\sinh d}{\sinh t}$.  Since
any triangle in hyperbolic space has area less than $\pi$, we have:
\[
\frac{d \vol(\delta_i)}{dt} < 4 \pi \frac{\sinh d}{\sinh t}
\]
Since $\sinh d$ is bounded, $\frac{d\vol(\delta_i)}{dt}$ decreases
exponentially with $t$. Applying the law of cosines to the triangle
with two sides of length $t$ and one side of length $i$ shows that
asymptotically $i$ is $2 t$ plus a constant, and so
$\frac{d\vol(\delta_i)}{di}$ decreases exponentially with $i$. \halmos
\end{proof}

By Lemma \ref{estimation}, if we fix $i_0$ then the set of $g \in X$
for which $\vol(\Str( f( g \cdot \sigma_i))) < \vol(\sigma_i) - i^2
\epsilon_i $ for some $i \geq i_0$ has measure less than $\mu(X)
\sum_{i = i_0}^{\infty} 1/i^2$.  Thus except on a small exceptional
set, we have 
\begin{equation}\label{vol_est}
\vol(\Str( f( g \cdot \sigma_i))) \geq \vol(\sigma_i) -
i^2 \epsilon_i
\end{equation}
for all $i \geq i_0$.

By Lemma \ref{hyperbolic}, $i^2 \epsilon_i \to 0$ as $i \to \infty$.
So letting $i_0 \to \infty$ we have that for almost all $g \in X$,
$\vol(\Str(f(g \cdot \sigma_i)))$ converges to $v_3$.  Let $r(i)$
denote the point of $r$ at a distance $i$ from $b$, the basepoint of
$r$.

  Let $\tau_i$ be a tetrahedron with vertices $b$, $r(i+1)$,
and the two vertices of $\sigma_i$ not on $r$ (see
Fig.~\ref{sig_n_tau}).  
\begin{figure}[h]
\centerline{\psfig{figure=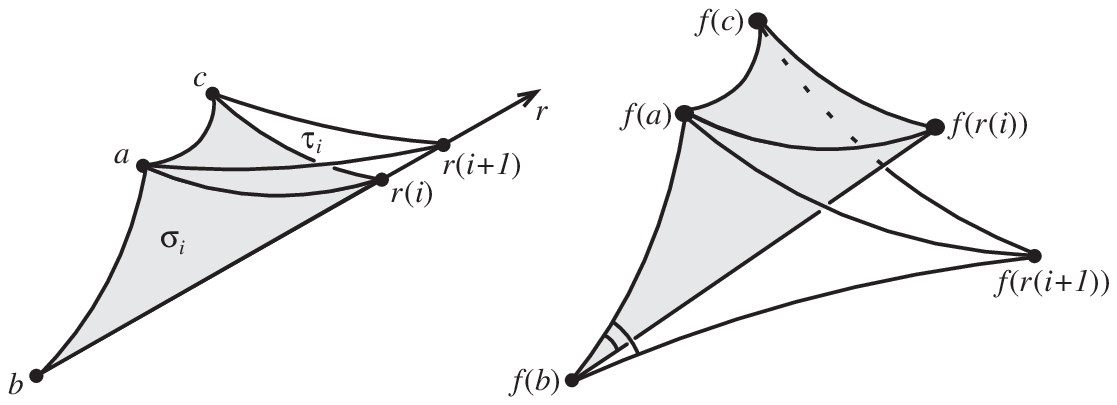}}
\firstcaption{The tetrahedra $\sigma_i$ and $\tau_i$ \label{sig_n_tau}}
\secondcaption{The tetrahedra $\Str(f(\sigma_i))$ and
  $\Str(f(\tau_i))$ \label{visual_angle}}
\end{figure}
Note that $v_3 - \vol(\sigma_i) > v_3 - \vol(\tau_i)$ and so the above
arguments show that for almost all $g \in X$, $\vol(\Str(f(g \cdot
\tau_i)))$ converges to $v_3$.  Suppose $g$ is such that both
$\vol(\Str(f(g \cdot \sigma_i)))$ and $\vol(\Str(f(g \cdot \tau_i)))$
converge to $v_3$.  For notational convenience, take $g$ to be the
identity.  I claim $f(r)$ converges to a point in $S_\infty^2$.  Since
$f$ is Lipschitz, it is enough to show that the $f(r(i))$ converge to
a point in $S_\infty^2$. Regular ideal tetrahedra are the only
tetrahedra of maximal volume, so since $\vol(\Str(f(\sigma_i)))$ goes
to $v_3$, we must have the distance from $f(b)$ to $f(r(i))$ going to
$\infty$ as $i$ goes to $\infty$.  Hence the $f(r(i))$ head out toward
$S_\infty^2$.  To show they converge, as opposed to wandering about
willy-nilly, we need to show the visual angle of $f(r(i))$ with
respect to $f(b)$ converges.  The change in visual angle between
$f(r(i))$ and $f(r(i+1))$ is the angle between the lines from $f(b)$
to $f(r(i))$ and from $f(b)$ and $f(r(i+1))$.  From
Fig.~\ref{visual_angle} we see this change is less than the sum of the
two indicated face angles of $\Str(f(\sigma_i))$ and
$\Str(f(\tau_i))$.

  The following lemma allows us to estimate these face angles:

\begin{Lemma}\label{angle} There is a constant $C > 0$ such that if
$\sigma$ is a tetrahedron with $\vol(\sigma)$ sufficiently close to
$v_3$, then for any face angle $\beta$ of $\sigma$:
\[
v_3 - \vol(\sigma) > C \beta^2
\]
\end{Lemma}

For large $i$, Eqn.~(\ref{vol_est}) shows $v_3 -
\vol(\Str(\sigma_i)) \leq (v_3 - \vol(\sigma_i)) + i^2 \epsilon_i =
\epsilon_i(1 + i^2)$.  From Lemma~\ref{angle} we have that the change in
visual angle for large $i$ is less than
 \[ 2\sqrt{\frac{\epsilon_i
     (1+i^2)}{C}}.
\]
By Lemma~\ref{hyperbolic} this is eventually exponentially decreasing
with $i$, and so the visual angles of the $f(r(i))$ converge.  Hence $f(g
\cdot r)$ converges for almost all $g \in X$.  Therefore, for almost
all geodesic rays $r$ in $\H^3$, $f(r)$ converges to a point in
$S_\infty^2$.

Moreover, as $f$ is Lipschitz any two rays which are asymptotic have
images under $f$ which converge to the same point of $S_\infty^2$.
Hence we have an extension of $f$, $\bar{f}:S_\infty^2 \to
S_\infty^2$, and you can check that $\bar{f}$ is measurable.  Modulo
the proof of Lemma~\ref{angle}, we have proven Claim~\ref{f_extends}.
 
Let us go back and prove the lemma.

\begin{proof}[Lemma \ref{angle}]
Let $v$ be the vertex of $\sigma$ which is the endpoint of the angle
$\beta$.  Without changing a neighborhood of $v$, push the other three
vertices of $\sigma$ to $S^2_\infty$ (this only decreases $v_3 -
\vol(\sigma)$).
\begin{figure}[th]
\centerline{ \psfig{figure=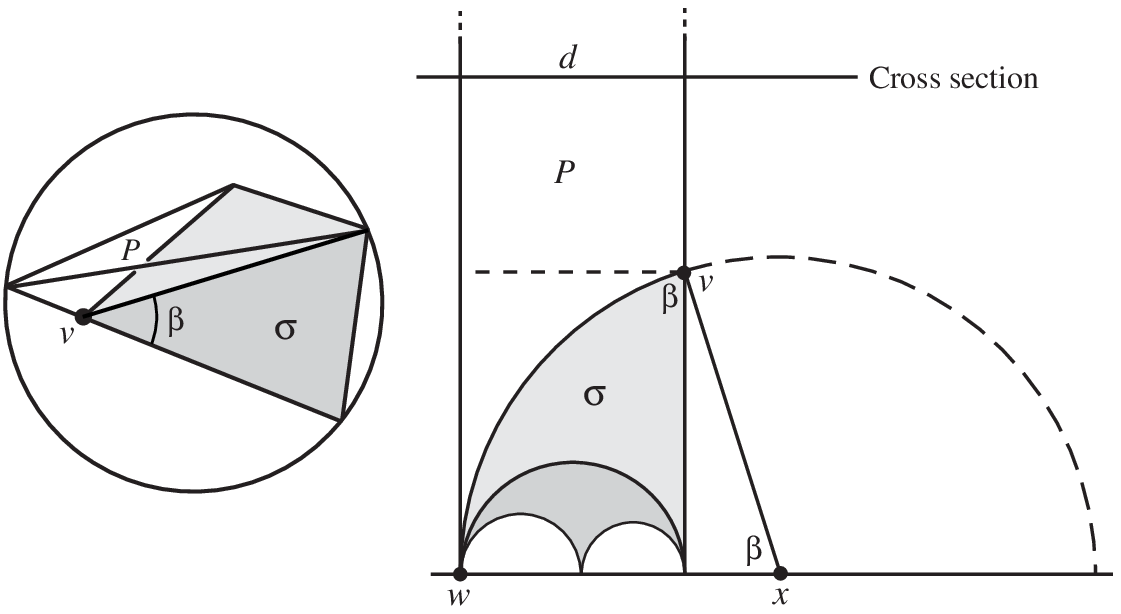}}
\firstcaption{Creating $P$ \label{create_P}}
\secondcaption{Estimating area of $P$\label{est_P}}
\end{figure} 
Extend an edge through $v$ which is a side of $\beta$ to $S^2_\infty$,
as in Fig.~\ref{create_P}.  Look at the part $P$ added on by doing
this.  Now $\vol(P) \leq v_3 - \vol(\sigma)$.  Consider
Fig.~\ref{est_P} in the upper half space model where we are using
Euclidean coordinates such that $\mbox{dist}(x,w) = \mbox{dist}(x,v) =
1$.

We will estimate the volume of $P$ above the dotted line.  By
requiring that $\vol(\sigma)$ be large, we can assume the dihedral
angles of $\sigma$ are close to $\pi/3$.  The indicated cross section
is then about an equilateral triangle whose area is bounded below by
$C_1 d^2$.  Note $d = 1 - \cos \beta$ with respect to our Euclidean
coordinate system, and so
\[
\vol(P) \geq \mbox{Vol.~above dotted line} = \int_{\sin \beta}^\infty
        \frac{ C_1 (1 - \cos \beta)^2 }{z^3} \, dz 
\]
Evaluating the integral we get
\[
     \vol(P)   \geq \frac{ C_1 (1 - \cos \beta)^2 }{2 \sin^2 \beta}
     \geq C_2 \beta^2
\]
for some constant $C_2$.  Thus $v_3 - \vol(\sigma) \geq   C_2 \beta^2$,
as desired.  \halmos
\end{proof}
 
Next, I will check Claim~\ref{pres_reg_tet} that $\bar{f}$ sends the
vertices of almost every positively oriented regular ideal tetrahedron
to the vertices of a positively oriented regular ideal tetrahedron.
Let $\delta_i$ denote a positively oriented regular tetrahedron with
side length $i$, with {\em center of mass} at a fixed point $b \in D$,
and whose vertices lie along four fixed geodesic rays $r_1, r_2, r_3,
r_4$ emanating from $b$.  Arguing as above, we can show that for
almost all $g \in X$, $\vol(\Str(f(g \cdot \delta_i)))$ converges to
$v_3$.  Moreover for almost all $g$ this is true and, in addition,
$f(g \cdot r_j)$ converges to a point $p_j \in S_\infty^2$ for all
$j$.  Since $\vol(\Str(f(g \cdot \delta_i)))$ converges to $v_3$, the
$p_j$ must span a regular ideal tetrahedron.  Since this is true for
almost all $g \in X$, the vertices of almost all regular ideal
tetrahedra are sent to regular ideal tetrahedra.  This proof of
Claim~\ref{pres_reg_tet}.

I will now prove Claim~\ref{done}, that $\bar{f}$ is essentially a
M\"obius transformation.  The space ${\cal T}$ of regular oriented
ideal tetrahedra with labeled vertices is a full measure subset of
$S^2_\infty \cross S^2_\infty \cross S^2_\infty$.  Let ${\cal T}^G$ be
the subset of ${\cal T}$ consisting of tetrahedra which $\bar{f}$
takes to regular oriented ideal tetrahedra.  We have just shown that
${\cal T}^G$ has full measure.  By Fubini's Theorem there is a $v_0
\in S^2_\infty$ such that almost all $T \in {\cal T}$ with first
vertex $v_0$ are in ${\cal T}^G$.  In fact, this is true for almost
all $v_0$, so we can assume that $\bar{f}(v_0)$ is defined (recall
that $\bar{f}$ is only defined by the process of looking at images of
geodesic rays for a full measure subset of $S_\infty^2$).

Without loss of generality, we can take both $v_0$ and $\bar{f}(v_0)$
to be the point at infinity in the upper half space model of $\H^3$.
Tetrahedra in ${\cal T}$ with first vertex at $\infty$ are equivalent
to oriented equilateral triangles in $\C$ with labeled vertices, which
are parameterized by a full measure subset of $\C \cross \C$.  It will
help the reader to think of $\C$ as having finite measure when
applying Fubini; we will only be concerned with which sets have
measure zero, a property which is invariant under diffeomorphism
\cite[Section~VI.1]{Boothby}.  For almost all lines $l$ through $0$,
almost all equilateral triangles with the edge between the first and
second vertices parallel to $l$ define tetrahedra which are in ${\cal
  T}^G$.  Assume that one such line is the real axis.  Let ${\cal S}$
denote tetrahedra with first vertex at $\infty$ and such that the edge
between the second and third vertices (the first and second vertices
of the corresponding triangle) is parallel to the real axis.

We know that ${\cal S}^G \equiv {\cal S} \cap {\cal T}^G$ has full
measure in ${\cal S}$.  Let $\omega$ be the $\sqrt[3]{-1}$ which has
positive imaginary part.  Then $\{ 0, 1, \omega \}$ is an oriented
equilateral triangle.  Let $L_0$ be all equilateral triangles in the
tiling of $\C$ by the triangle $\{ 0, 1, \omega \}$.  Let $L_k$ be the
same set of triangles scaled by $2^{-k}$.  Let $L = \bigcup_{k \in \Z}
L_k$ be this nested family of equitriangular lattices (See
Fig.~\ref{lattice}).
\begin{figure}[ht]
\centerline{ \psfig{figure=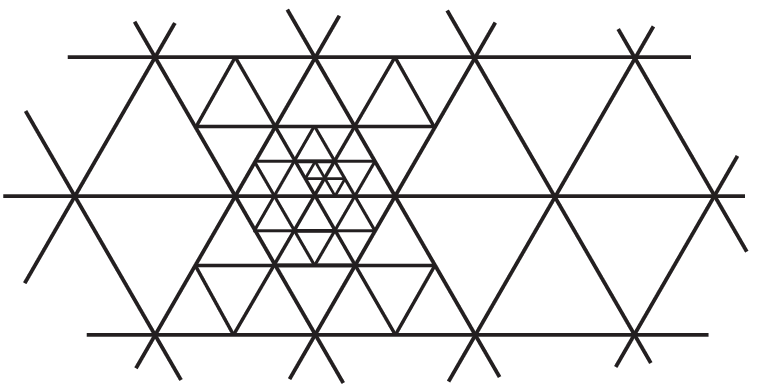}}
\caption{Some of the triangles in the nested family of lattices $L$}
\label{lattice}
\end{figure} 

I claim there is an $r \in \R$ such that for almost all $z \in \C$,
the entire countable set of triangles $z + rL$ are in ${\cal T}^G$.
Consider the submersion $\pi \maps \C \cross \R \cross \Z \cross \Z
\cross \Z \to {\cal S}$ which sends $(z, r, k, n, m)$ to the
equilateral triangle with vertices
\[ 
\left(z + r 2^{-k} (n + m \omega), z + r 2^{-k} (n + 1 +
m \omega), z + r 2^{-k}(n + (m + 1) \omega)\right)
\]
in $z + r L_k$.  We will think of $\Z$ as having a measure where the
measure of $q$ is $1/q^2$.  As $\pi$ is a submersion, $\pi^{-1}({\cal
  S}^G)$ has full measure.  Thus by Fubini, for almost all $r$ and
$z$, we have $\pi^{-1}\left({\cal S}^G\right) \cap \left(\{r \} \cross \{z\} \cross \Z
\cross \Z \cross \Z\right)$ has full measure, that is equal to $\{r \}
\cross \{z\} \cross\Z \cross \Z\cross \Z$, as desired.  Without loss
of generality assume $r = 1$ has this property.  So for almost all $z
\in \C$ all triangles in $z + L$ are in ${\cal S}^G$.  This forces
$\bar{f}(z + L)$ to be family of nested equitriangular lattices (see
Fig.~\ref{lattice}).  For each $z$ there is a complex number $h(z)$
such that:
\[
\bar{f}\left(z + 2^{-k}( n + m \omega)\right) = \bar{f}(z) + h(z)
2^{-k}
(n + m \omega)
\]
for all $\{n, m, k \} \subset \Z$.  I claim the function $h$ is
invariant under the group of translations of the form $z \mapsto z +
2^{-j}(a + b \omega)$, where $\{j, a, b\} \subset \Z$.  Let $z' = z +
2^{-j}(a + b \omega)$. We have
\begin{eqnarray}
\bar{f}(z' ) &=& \bar{f}(z) + h(z) 2^{-j}( a +  b \omega) 
\label{first} \\
\bar{f}(z') + h(z')  2^{-j} &=& \bar{f}(z) + h(z)  2^{-j}( a + 1 + b\omega) =
 \bar{f}(z' + 2^{-j}) \label{second}
\end{eqnarray}
Subtracting Eqn.~(\ref{first}) from Eqn.~(\ref{second}) we get $h(z') =
h(z)$.  Our group of translations is dense, and so acts ergodically.
Therefore $h$ is constant almost everywhere.  But then $\bar{f}(z') =
\bar{f}(z) + h  2^{-j}(a + b \omega)$ almost everywhere which implies
that $\bar{f}(z) - h \cdot z$ is invariant under our group of
translations.  So there is a constant $c$ such that $\bar{f}(z) - h z =
c$
almost everywhere and thus $\bar{f}(z) = c + h z$ almost
everywhere.

 Thus $\bar{f}$ is essentially a M\"obius transformation.  The
corresponding
M\"obius transformation conjugates the actions of $\rho_0$ and
$\rho_1$ on $S^2_\infty$ to one another.  Thus $\rho_0$ and $\rho_1$
are conjugate, and $\rho_1$ is discrete and faithful. \halmos
\end{proof}

\bibliographystyle{math}

\bibliography{standard}

\end{document}